\documentclass{article}
\usepackage{enumerate}
\usepackage{myMac}
\draftdim

\newcommand{\calcc}{\ol{\calC}}

\newcommand{\calc}{\calC}

\newcommand{\caln}{\calN}
\newcommand{\calv}{\calV}
\newcommand{\calno}{\calN(I_0)}
\newcommand{\calmo}{\calM_0}
\newcommand{\cali}{\calI}

\newcommand{\caltp}{\calT_{P}}
\newcommand{\tcc}{T_\calc}
\newcommand{\calto}{T_{I_0}}
\newcommand{\calco}{\calc_{0}}

\newcommand{\fk}{f^{(k-1)}}
\newcommand{\annc}{\calA_\calc}
\newcommand{\calic}{\calI_\calc}
\newcommand{\calip}{\calI_p}
\newcommand{\nie}{\texttt{Newton-Incl-Exc}~}
\newcommand{\nisol}{\texttt{Newton-Isol}~}

\newcommand{\gxdx}{G(x) dx}
\begin{document}
\title{Near Optimal Subdivision Algorithms for Real Root Isolation}
\author{
%
% 1st. author
Vikram Sharma\\
{Institute of Mathematical Sciences}\\
{Chennai, India 600113}\\
{vikram@imsc.res.in}\\
\and
% 2nd. author
Prashant Batra\\
{TU Hamburg}\\
{Hamburg, Germany}\\
{batra@tuhh.de}
}%author

\maketitle
\begin{abstract}
The problem of isolating real roots
of a square-free polynomial inside a given interval $I_0$ is a fundamental problem.
Subdivision based algorithms are a standard approach to solve this problem.
Given an interval $I$, such algorithms rely on two predicates: an exclusion predicate, 
which if true means $I$ has no roots, and an inclusion predicate, which if true, reports an isolated root in $I$.
If neither predicate holds, then we subdivide the interval and proceed recursively, starting from $I_0$.
Example algorithms are Sturm's method (predicates based on Sturm sequences), 
the Descartes method (using Descartes's rule of signs),
and \texttt{Eval} (using interval-arithmetic). For the canonical problem of isolating all
real roots of a degree $n$ polynomial with integer coefficients of bit-length $L$, the subdivision tree
size of (almost all) these algorithms is bounded by $O(n(L+\log n))$. This is known to be optimal
for subdivision based algorithms.

We describe a subroutine that improves the running time of any subdivision algorithm
for real root isolation. The subroutine first detects clusters of roots using a result of Ostrowski,
and then uses Newton iteration to converge to them.
Near a cluster, we switch to subdivision, and proceed recursively. 
The subroutine has the advantage that it is independent
of the predicates used to terminate the subdivision.
This gives us an alternative and simpler approach to recent developments of Sagraloff (2012)
and Sagraloff-Mehlhorn (2013), assuming exact arithmetic.

The subdivision tree size of our algorithm
using predicates based on Descartes's rule of signs is bounded by $O(n\log n)$,
which is better by $O(n\log L)$ compared to known results. 
Our analysis differs in two key aspects. First, 
we use the general technique of continuous amortization 
from Burr-Krahmer-Yap (2009), and 
second, we use the geometry of clusters of roots
instead of the Davenport-Mahler bound. The analysis
naturally extends to other predicates.

\end{abstract}

%A category including the fourth, optional field follows...
%\category{D.2.8}{Software Engineering}{Metrics}[complexity measures, performance measures]


%\keywords{Real root isolation, Subdivision algorithms, Newton diagram, Continuous amortization, Integral analysis.}

%%%%%%%%%%%%%%%%%%%%%%%%%%%%%%%%%%%%%%%%%%%%%%%%%%
\section{Introduction}
%%%%%%%%%%%%%%%%%%%%%%%%%%%%%%%%%%%%%%%%%%%%%%%%%%
Given a polynomial $f \in \RR[x]$ of degree $n$, the problem is to isolate the
real roots of $f$ in an input interval $I_0$, i.e., 
compute disjoint intervals which contain exactly one real root of 
$f$, and together contain all roots of $f$ in $I_0\cap \RR$. 
% There are many variants and known
% algorithms for the problem.
Subdivision based algorithms
have been successful in addressing the problem. 
A general subdivision algorithm uses two predicates,
given an interval $I$: the exclusion predicate $C_0(I)$, which if true means
 $I$ has no roots; the inclusion predicate, $C_1(I)$, which if true means $I$ has exactly one root. 
The algorithm outputs a \dt{root-partition} $\calP$ of $I_0$, i.e., a set of pairwise disjoint open intervals
such that for each interval either $C_0$ or $C_1$ holds, and $I_0 \sm \calP$ contains no roots of $f$.
To compute isolating intervals for roots of $f$, check the sign of $f$ at the endpoints of the intervals
in $\calP$ (this works if $f$ is square-free). The following generic subdivision algorithm constructs
a root-partition:
\progb{
\texttt{Isolate}$(I_0)$\\
{\sc Input:} $f \in \RR[x]$ and an interval $I_0 \ib \RR$.\\
{\sc Output:} A root-partition $\calP$ of $f$ in $I_0$.\\
0.\> Preprocessing step.\\
1.\> Initialize a queue $Q$ with $I_0$, and $\calP \ass \es$.\\
2.\> While $Q$ is not empty\\
\>\> Remove an interval $I=(a,b)$ from $Q$.\\
\>\> If $C_0(I) \vee C_1(I)$ then add $I$ to $\calP$.\\
\>\> else\Comment{Subdivide $I$}\\
\>\>\> Let $m \ass (a+b)/2$.\\
\>\>\> Push $(a, m)$ and $(m,b)$ into $Q$.\\
3.\> Output $\calP$.
}
The algorithm is guaranteed to terminate for square-free polynomials; 
otherwise we get an infinite sequence of intervals converging to a root of multiplicity
greater than one.
Some standard choices of the predicates and the corresponding algorithms are:
\begin{tightenum}{r}
\item  Sturm sequences and Sturm's method \cite{davenport:85},
\item Descartes's rule of signs and the Descartes method \cite{collins-akritas:76},
\item  Interval-arithmetic based approaches and \texttt{Eval} \cite{burr-krahmer-yap:continuousAmort:09}.
\end{tightenum}
The complexity of these algorithms is well understood for the
benchmark problem of isolating all real roots of a square-free integer polynomial with coefficients
of bit-length $L$. One measure of complexity is the size of the subdivision tree constructed by
the algorithm. For the first two algorithms a bound of $O(n(L+\log n))$ was shown in
\cite{davenport:85} and \cite{eigenwillig-sharma-yap:descartes:06}, respectively. 
For \texttt{Eval} a weaker bound of $O(n(L+n))$ was established in \cite{sharma-yap:near-optimal:12}.
It is also known \cite{eigenwillig-sharma-yap:descartes:06} that the bound $O(n(L+\log n))$
is essentially tight for any algorithm doing uniform subdivision, i.e., 
reduces the width at every step by some constant (in our case by half).
%\footnote{Here we are talking of uniform
%subdivision, i.e., the interval is subdivided into equal halves. Examples of a non-uniform subdivision 
%algorithm are the continued fraction based approaches, but even for them similar lower bounds
%have been derived \cite{collins-krandick:12}.} 

Uniform subdivision cannot improve on the bound mentioned above because
it only gives linear convergence to a ``root cluster'', i.e., roots which
are relatively closer to each other than to any other root. But it is known that
from points sufficiently far away from the cluster, 
Newton iteration (more precisely, its variants for multiple roots) 
converges quadratically to the cluster. 
This has been an underlying idea in improving the linear convergence
 of subdivision algorithms for root isolation \cite{pan:approx-poly-zeros:00,sagraloff-newton+descartes:12,sagraloff-mehlhorn:13},
and has also been combined with homotopy based approaches \cite{yakoubsohn:zero-clusters:00,sasaki-terui:cluster:09}.

We follow the same idea with some key differences.
Given $C_0$ and $C_1$,
our algorithm can be described as follows (we only 
give the inner loop, see \refSec{algorithm} for complete details):
\progb{
\texttt{Newton-Isol}$(I_0)$\\
% {\sc Input:} $f(x) \in \RR[x]$ and an interval $I_0 \ib \RR$.\\
% {\sc Output:} A root-partition $P$ of $f$ in $I_0$.\\
% 0.\> Preprocessing step.\\
% 1.\> Initialize a queue $Q$ with $I_0$, and $P \ass \es$.\\
% 2.\> While $Q$ is not empty\\
% \>\> Remove an interval $I=(a,b)$ from $Q$.\\
\> \dots\\
\> If $C_0(I)\vee C_1(I)$ then add $I$ to $\calP$.\\
\> else if a cluster $\calc$ of roots is detected in $I$ then\\
\>\> Apply Newton iteration to approximate $\calc$\\
\>\> while quadratic convergence holds.\\
\>\> Estimate an interval $J$ containing $\calc$.\\
\>\> Push $J$ into $Q$.\\
\> else \Comment{Subdivide $I$}\\
% \>\> Let $m \ass (a+b)/2$.\\
% \>\> Push $(a, m)$ and $(m,b)$ into $Q$.\\
\>\> \dots
% 3.\> Output $P$.
}
% (i) If $C_0(I)$ or $C_1(I)$ holds then discard or output $I$.\\ 
% (ii) If a cluster of roots is detected in $I$ then \\
% (iii) If all the tests above  fail, then subdivide $I$ into equal halves and proceed recursively.\\
% Note that the second step is the only key difference from the generic framework of a subdivision
% algorithm given above. 
For detecting root clusters, we use a result of Ostrowski
based on Newton diagram of a polynomial  \cite{ostrowski:graeffe1:40}; other choices are
based on a generalization of Smale's $\alpha$-theory (see \cite{giusti+2:zeros-analytic:05}
and the references therein); the details can be found in \refSec{notation}.
These tools and approaches have been used earlier
\cite{pan:approx-poly-zeros:00}, however, our approach has the following differences:
\begin{tightenum}{r}%[label={(\roman*)},itemsep=-10pt,topsep=0pt,parsep=0pt]
\item The tools used to detect and estimate the size of a cluster are independent
of the particular choice of the exclusion-inclusion predicates 
(cf. \cite{sagraloff-newton+descartes:12}). 
This way we obtain a general approach to improve any subdivision algorithm.
\item Another difference is the method that is combined with bisection to improve convergence.
In \cite{sagraloff-newton+descartes:12} Abbott's QIR method 
is combined with the Schr\"oder operator \cite{giusti+2:zeros-analytic:05}, whereas  we
apply standard Newton iteration to a suitable derivative of $f$. 
The former combination is a backtracking approach to get quadratic convergence; the latter gives quadratic
convergence right away (but perhaps increasing subdivisions). This has the
advantage of separating the Newton iteration steps from the subdivision tree,
which is reflected in the bounds on the subdivision tree size for the two approaches: 
for the former we have $O(n\log (nL))$, and for the latter we have $O(n\log n)$.
The number of quadratically converging steps remains the same in both cases.
\item  Our approach can be modified to
isolate complex roots; replace binary subdivision with a quad-tree subdivision, and choose
appropriate predicates (e.g., Ostrowski's result mentioned above, or Pellet's test). This avoids
Graeffe iteration (cf. \cite{pan:approx-poly-zeros:00}),
and yet the modified algorithm can be shown to attain a near optimal bound on subdivision tree size.
\end{tightenum}
% - The cluster detection result of Ostrowski works in the more general setting of analytic functions.
% Thus if  we choose interval arithmetic based predicates, e.g., those suggested in \cite{sagraloff-sharma-yap:analytic:13},
% then we obtain an algorithm for isolating roots of analytic functions.
% One way to bound the bit-complexity of \texttt{Newton-Isol} is to bounds its
% algebraic complexity and precision requirements.
%In this paper, we assume the Real RAM model, and bound the 
%algebraic complexity of \texttt{Newton-Isol}. This involves two steps:
%first derive an upper bound on the
%size of the subdivision tree; second, bound the worst case number of operations at each 
%node in the tree; multiply the two bounds to obtain an upper bound on the algebraic cost.

In this paper, we focus on bounding the size of subdivision
tree of \texttt{Newton-Isol}.
% assume the Real RAM model and bound the arithmetic complexity of 
% \texttt{Newton-Isol}. An important component of the analysis is to bound the size of subdivision
% tree constructed by the algorithm; the bound on arithmetic complexity follows by
% multiplying the size with the worst case complexity of a node in the tree.
For this purpose, we  use the general framework
of continuous amortization \cite{burr-krahmer-yap:continuousAmort:09,burr:contamort:13}.
%which provides a unified approach to study the complexity of subdivision algorithms.
The key idea here is to bound the tree size
by $\int_{I_0} G(x) dx$, where $G$ is a suitable ``charging'' function 
corresponding to the predicates used in the algorithm (e.g., see \cite{burr:contamort:13}).
Our key contributions are the following:
\begin{tightenum}{r}
\item We derive a near optimal bound of $O(n\log n)$ on the size of the subdivision tree 
of \texttt{Newton-Isol} when $C_0$, $C_1$ are based on Descartes's rule of signs
(see \refThm{bound}). This is the first application of the continuous amortization framework
to a non-uniform subdivision algorithm.
\item 
We show that if the distance of the cluster center to the nearest root outside the cluster
exceeds roughly $n^3$ times the diameter of the cluster,
then Ostrowski's criterion for cluster detection works,
and we obtain quadratic convergence to the cluster center (see \refLem{converse}).
\item  Our analysis crucially uses the cluster tree of the polynomial (see \refPro{ct}).
We derive an integral bound on the size of the subdivision tree (see \refThm{intbound}).
The usual approach to upper bound this integral is to break it over
the (real) Voronoi regions of the roots \cite{burr:contamort:13}.
We instead break the integral over the Voronoi regions corresponding to the clusters 
in an inductive manner based on the cluster tree.
The integral over the portion of the region outside the cluster is bounded using known techniques.
However, for the portion inside the cluster, we devise an amortized bound on the integral (see \refLem{dense}), which is
of independent interest, and is analogous to the improvement given by Davenport-Mahler bound
over repeated applications of the root separation bound. It is this result that underlies the $O(n\log n)$ bound.
A simple argument extends these bounds to Sturm's method and the \texttt{Eval} algorithm.
The details  are in \refSec{complexity}. 
\end{tightenum}
% Since the
% cluster tree has at most $O(n)$ nodes, we obtain the bound of $O(n \log n)$ on the size
% of the subdivision tree. We show our analysis using the Descartes's rule of signs 
% as the predicate of choice, though the argument can be extended for other predicates
% as well.

\section{Notation and Basic Results}\label{sec:notation}
Let $f \in \RR[x]$  be a square-free polynomial of degree $n\ge 2$ 
and $Z(f)\ibp \CC$ be its set of roots. 
Given a finite pointset $S\ib \RR^2$, let $D_S$ be the disc $D(m_S,r_S)$
such that $m_S$ is the centroid of the points in $S$,
and $r_S$ is the least radius such that all the points in $S$ are contained in $D(m_S, r_S)$.
Given a $\lambda \in \RR_{> 0}$, define $\lambda D_S \as D(m_S, \lambda r_S)$.
We borrow the following definition from \cite{sagraloff-sharma-yap:analytic:13}: 
A subset $\calc \ib Z(f)$ of size at least two  is called a (root) \dt{cluster}  if 
the only roots in $3D_\calc$ are from $\calc$.
%\footnote{An alternative definition could be that the maximum distance 
%between any pair of points in $\calc$
%is much smaller than the minimum distance between any point in $\calc$ and $Z(f)\sm \calc$.}
We treat individual roots as (trivial) clusters. In this paper,
the non-real roots in $\calc$ come in conjugate pairs.
Therefore, the center of $D_\calc$ will always be in $\RR$.
Define $R_\calc$ as the distance from $m_\calc$ to the nearest point in the set $Z(f)\sm \calc$.
From the definition it follows that $Z(f)$ trivially forms a cluster and $R_{Z(f)}=\infty$. 
Given an interval $I$, let $m(I)$ denote its midpoint and $w(I)$ its width.
We will often use the shorthand $I=[m(I)\pm w(I)/2]$, and for $\lambda > 0$,
$\lambda I \as [m(I)\pm \lambda w(I)/2]$.
An \dt{interval $I$ contains a cluster} $\calc$ if $\calc \ib D(m(I), w(I)/2)$.

We use the following convenient notation in the subsequent definitions: 
for $x, y \in \RR$, `$x \gg y$' if there is a constant $c \ge 1$
such that $x \ge c y$; similarly define $x \ll y$.

% The \dt{spread} $\sigma_\calc$ 
% of a cluster is defined to be the ratio $R_\calc/r_\calc$.
A \dt{strongly-separated cluster (ssc)} is a cluster $\calc$ for which 
$R_\calc/r_\calc \gg n^3$;  the exact
constant can be found in \refCor{con}. For a ssc $\calc$, define the following three quantities:
\begin{tightenum}{r}
\item The interval $I_\calc\as [m_\calc\pm c\cdot kr_\calc]$, for some constant $c\ge 1$.
\item The interval $\calic \as \set{x: |x-m_\calc|\ll R_\calc/n^2}$.
\item The annulus $\annc \as \calic \sm I_\calc=\set{z \in \CC: |\calc|r_\calc \ll |z-m_\calc|\ll R_\calc/n^2}$.
\end{tightenum}
The exact constants in these definitions are given in \refLem{converse}.
See \refFig{ssc} for an illustration of these concepts.
If $\calc$ is not a ssc, then we define 
$ I_\calc \as [m_\calc \pm r_\calc]$ and $\calic\as 2I_\calc$. 
Note that for all clusters $\calc$,  $I_\calc \ib \calic$.
We will need the following result later in our analysis \cite[Lemma 2.1]{sagraloff-sharma-yap:analytic:13}:
\bprol{ct}
Given a root cluster $\calc$ of $f$. There is a unique unordered tree $T_\calc$
rooted at $\calc$ whose set of nodes are the clusters contained in $\calc$, and
the parent-child relation is subset inclusion. Let $T_f$ be the tree where
the parent is the cluster $Z(f)$ of all roots.
\eprol
The result originally is stated for root clusters of $f \in \CC[x]$. However, for $f \in \RR[x]$
the clusters come in conjugate pairs, and by taking the union of such pairs the result
still holds. The tree $T_\calc$ is called the \dt{cluster tree of $\calc$}. The leaves of this tree are the roots
in $\calc$.

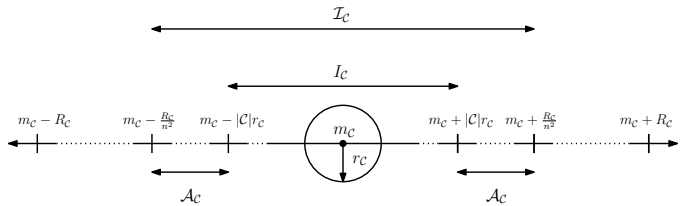
\vfigpdf{Geometry of a ssc $\calc$. We focus on the the relative geometry,
overlooking the exact constants involved in the definition of the intervals. 
% The result in \cite{pawlowski:zeros-of-derivatives:99} roughly states that
% for $j < |\calc|$, $f^{(j)}$ has $|\calc|-j$ roots in $D(m_\calc, r_\calc)$
% and the remaining $n-|\calc|$ roots outside $D(m_\calc, R_\calc/(4n^2))$.
}{ssc}{0.45}

% A \dt{maximal cluster} in $T_\calc$ is a (non-trivial) cluster  
% at depth one in $T_\calc$. A \dt{minimal cluster} in $T_\calc$ is a cluster all of whose children
% are leaves. 

% Let $B(I)$ denote the box centered at $m(I)$ with width $w(I)$.

%%%%%%%%%%%%%%%%%%%%%%%%%%%%%%%%%%%%%%%%%%%%%%%%%%
\subsection{Cluster Detection and Approximation}
The literature on detection and approximation of root clusters is vast
(see \cite{giusti+2:zeros-analytic:05} and the references therein).
One approach is based on Pellet's test: if for a complex polynomial $f(x)=\sum_{i=0}^na_i x^i$
there is an $r > 0$ such that $|a_k|r^k > \sum_{i \neq k} |a_i|r^i$ then the disc
$D(0, r)$ contains exactly $k$ roots of $f$. A point $z \in \CC$
is said to satisfy Pellet's test, if there is a $k$ and $r$
for which the test holds with the coefficients of $f(x+z)$.
Results in \cite{yakoubsohn:zero-clusters:00,giusti+2:zeros-analytic:05} 
generalize Smale's $\alpha$-theory and relate it to Pellet's test; 
an alternative derivation based on tropical algebra is given in \cite{sharify:thesis}.
We instead use a result by Ostrowski \cite{ostrowski:graeffe1:40}.
% one reason being the computational effectiveness of the test.

We need the following definitions.
Let $f(x) = \sum_{i=0}^n a_i x^i$, where $a_i \in \CC$.
%for simplicity, assume that $a_0 a_n \neq 0$. 
With each index $i$, $a_i \neq 0$, associate the point
$P_i \as (i, -\log |a_i|) \in \RR^2$. The lower-hull of the convex-hull of these points
is called the \dt{Newton diagram} of $f$. % note that $P_0$ and $P_n$  belong to the diagram.
Given an index $k\in \set{0\dd n}$, let $y_k$ be the point such that $(k, y_k)$ is on the diagram. 
%Define $\tau_k \as e^{-y_k}$, for $k \in \set{0 \dd n}$. %, and $0$ otherwise.
Define $\rho_k \as e^{y_k-y_{k-1}}$, for $1 \le k\le  n$, %the  $k$th \dt{exponential slope}
$\rho_{n+1}\as \infty$, and the 
$k$th \dt{deviation} $\Delta_k \as \rho_{k+1}/\rho_k$, for $0 < k < n$.

Let $\alpha_1 \dd \alpha_n \in \CC$ be the roots of $f$ ordered such that 
$|\alpha_1| \le |\alpha_2| \le \cdots \le |\alpha_n|$. Ostrowski showed the following
fundamental relation between the absolute values of the roots and $\rho_k$'s
\cite[p.~143]{ostrowski:graeffe1:40}:
        \beql{ost}
        \frac{1}{2k} < \frac{|\alpha_k|}{\rho_k}< 2(n-k+1).
        \eeql
Given $z\in \CC$, we will be interested in the Newton diagram of $f(x+z)$.
If $f_j(z) \as f^{(j)}(z)/j!$, then from a result of
Ostrowski \cite[p.~128]{ostrowski:graeffe1:40} we get:
% Ostrowski also gave a closed form expression for the exponential slopes
% \cite[p.~128]{ostrowski:graeffe1:40} at two true corner points $P_k,P_{k+1}$
% of the Newton diagram:
        \beql{rkr}
        \rho_k(z) = \max_{j<k} \abs{\frac{f_j(z)}{f_k(z)}}^{\frac{1}{(k-j)}},\text{ and }
        \rho_{k+1}(z) = \min_{j>k} \abs{\frac{f_k(z)}{f_j(z)}}^{\frac{1}{(j-k)}}.
        \eeql
The RHS is defined for any $k$ such that $f_k(z) \neq 0$; however, we 
are only interested in those $k$ for which $P_k$ is 
on the diagram. The $k$th deviation $\Delta_k(z)\as \rho_{k+1}(z)/\rho_{k}(z)$.
We have the following result for detecting clusters:
\bleml{cluster}
If $\Delta_k(z) \ge 27$, for some index $0 < k< n$, then 
there are exactly $k$ roots in $D(z, 3\rho_k(z))$ and $D(z, \rho_{k+1}(z)/3)$. Moreover,
as $\rho_{k+1}(z)/3 \ge 9\rho_k(z)$, these roots form a cluster.
\eleml
The proof shows that the inequality $\Delta_{k}(z)\ge 27$ implies that Pellet's
test holds for $D(z, r)$, $3\rho_k(z) \le r \le D(z, \rho_{k+1}(z)/3)$ 
(see \cite[Thm.~1.5]{giusti+2:zeros-analytic:05}). 
Since the $P_i$'s are sorted by x-coordinate, 
all the $\rho_k$'s can be computed in $O(n)$ steps using, e.g., Graham's scan for convex hull computation. 
% Given $z \in \CC$,
% let $\rho_k(z)$ and $\Delta_k(z)$ be defined as above but with respect to the Newton diagram
% of the shifted polynomial $f(x+z)$.
% The result above gives us a test to detect root clusters near $z$ using $\Delta_k(z)$.

Once we have detected a cluster $\calc$ near $z$, we want a good approximation to 
$m_\calc$. A standard way is to do the iteration $z_{i+1}= z_i - kf(z_i)/f'(z_i)$, starting
from $z$, 
but this may not be numerically desirable, as both $f$ and $f'$ are small near 
$\calc$. Another option is to use the standard Newton iteration applied to $\fk$.
We show that if $\Delta_k(z) \ge 27$, then $z$ is an approximate zero,
in the sense of Smale et al.~\cite[p.~160, Thm.~2]{bcss:bk},  
to the root of $\fk$  in $D(z, \frac{3\rho_k(z)}{2k})$.
% We recall the three
% basic functions from there: $\beta(g, z) \as  \abs{\frac{g(z)}{g'(z)}}$, 
% $\gamma(g,z) \as \max_{j >1}\abs{\frac{g^{(j)}(z)}{j! g'(z)}}$, and $\alpha(g,z) =\beta(g,z)\gamma(g,z)$.
% Let $\beta_k(z), \gamma_k(z), \alpha_k(z)$ be the  corresponding quantities for $g\as \fk$.
% We have the following relation:
% \beql{alphakrho}
% \alpha_k(z) \Delta_k(z) \le 2.
% \eeql
% This means that if the $k$th deviation is sufficiently large, then $z$ is an approximate zero to
% the root of $\fk$ in the disc $D(z, 1.5\beta_k(z))$ \cite[p.~160, Thm.~2]{bcss:bk}.
Subsequently we show that if $\Delta_k(z)\ge c_0$, for some constant $c_0\ge 27$, 
then for all $z'$ in this disc $\Delta_k(z') \ge 27$, and hence there is a cluster of $k$ roots in $D(z', 3\rho_k(z'))$.
Moreover, the cluster is exactly $\calc$. %  since
% $3 \rho_k(z') \le 6e^6\rho_k(z) < \rho_{k+1}(z)/3.$
% \bleml{deltak}
% If $\Delta_{k}(z) > 16$, then for all $z' \in D(z, 1.5\rho_k(z))$ we have 
% $\Delta_k(z')\ge \Delta_k(z)/6e^6$.
% \eleml
% Thus then for all 
% $z' \in D(z, 1.5\beta_k(z))$, $\Delta_k(z') \ge 27$
% where in the last step we use the fact that $\Delta_k(z) > 2e^6\times 9$; 
% Instead of explicitly mentioning the
% constant, we will use the notation $\Delta_k(z) \gg 1$. %  to mean that there is a 
% universal constant $c_0$ such that $\Delta_k(z) \geq c_0$.
These results are summarized in the following:
%\vspace*{-.5cm}
\bleml{correctness}
Let $z \in \CC$ be such that $\Delta_k(z) \ge c_0$, for some $k \ge 2$, 
$\calc$ be the cluster in $D(z, 3\rho_k(z))$, and $D' \as D(z, \frac{3\rho_k(z)}{2k})$.
Then the following hold:
%(i) There is a cluster $\calc$ of $k$ roots in $D(z, 3\rho_k(z)) \ib D(z, \frac{\rho_{k+1}(z)}{3})$.\\
\begin{tightenum}{r}
\item  $z$ is an approximate zero to the root $z^*$ of $\fk$ in $D'$ and
 the Newton iterates starting from $z$ are in $D'$.
\item  For all $z' \in D'$,
%$\rho_k(z') \ll \rho_k(z)$ and $\rho_{k+1}(z')\gg \rho_{k+1}(z)$, and hence 
$\Delta_k(z') \ge 27$, and % As $D(z', 3\rho_k(z')) \ibp D(z, \frac{\rho_{k+1}(z)}{3})$, 
$\calc$ is the cluster in $D(z', 3\rho_k(z'))$.
\item  If $z, w$ are such that $\Delta_k(z), \Delta_k(w) \ge 27$ and $D(z, 3\rho_k(z))$,
$D(w, 3\rho_k(w))$ intersect then the discs have the same cluster.
\end{tightenum}
\eleml
%\vspace*{-.5cm}
The proof is given in the appendix. 
We choose $c_0 \as 27\times 6e^6$.
Given $z\in \CC$, a value of $k$ satisfying the condition $\Delta_k(z) \ge c_0$ is called an \dt{admissible
value} for $z$, with the corresponding \dt{inclusion disc} $D(z, 3\rho_k(z))$.
% that contains a cluster of size $k$.
Note that there can be more than one admissible value for a point $z$ corresponding to
clusters of different sizes.

\section{The Algorithm}\label{sec:algorithm}
Let $C_0$ and $C_1$ be some exclusion and inclusion predicate respectively.
The following algorithm takes as input $f$ and an interval $I_0$ and outputs a root partition
of $I_0$. 
%\vspace*{-1cm}
\progb{
\texttt{Newton-Isol}$(f, I_0)$\\
1\> Initialize  $\calP \ass \es$, $\Phi \ass \es$; let $Q$ be an empty queue.\\
1.a.\> If this is a recursive call then subdivide $I_0$ and \\
\> push the two halves into $Q$;  else $Q \ass \set{I_0}$.\\
%\> \Comment{$R$ contains intervals to be considered recursively}\\
2.\> While $Q$ is not empty do\\
\>\> Remove an interval $I$ from $Q$.\\
2.a.\>\> If $C_0(I)\vee C_1(I)$ then add $I$ to $\calP$.\\
\>\>else if \texttt{Newton-Incl-Exc}($I$) is successful then\\
\>\>\> Let $(J,k)$ be the pair returned.\\
2.b.\>\>\> If $\forall J' \in \Phi$, $J\si J'=\es$ and $J \si I_0 \neq \es$ then\\
%\>\>\> \Comment{I.e., this cluster has not occurred earlier and overlaps the input interval}\\
2.c.\>\>\>\>  $\forall$ $I'\in Q$, $I'\ass I'\sm D(m(J),\frac{\rho_{k+1}(m(J))}{3})$.\\
2.d.\>\>\>\> Add $J\si I_0$ to $\Phi$.\\
%\>\> else \Comment{i.e., \texttt{Newton-Incl-Exc} returned failure}\\
\>\> else subdivide $I$ and push the two halves into $Q$.\\
3.\> Return $\calP \su_{J \in \Phi} \texttt{Newton-Isol}(f, J)$.
}
%\vspace*{-.5cm}
The input to \texttt{Newton-Incl-Exc} is an interval $I=(a,b)$. If the predicate is successful then
it returns  an interval $J$ containing a cluster such that $w(J) < w(I)/2$, and an 
admissible value $k$ for $m(J)$; otherwise it returns failure.
\progb{
\texttt{Newton-Incl-Exc} $(f,I)$\\
1.\> Let $m \as (a+b)/2$.\\
2.\> For $p \in \set{a,m,b}$, let $k_p \ge 2$ be the {\em smallest} admissible\\
\> value $k$ for $p$ such that $I\ib D\paren{p, \frac{\rho_{k+1}(p)}{3}}$. \\
%\> \Comment{This is done by computing respective Taylor coefficients and a Newton diagram computation.}\\
%.\> Let $P\ib \set{a,m,b}$ be the set of points for which $k_p<n$.\\
3.\> If the three admissible values are equal and the three\\
\> inclusion discs are contained in $D\paren{m, \frac{\rho_{k_m+1}(m)}{3}}$ then:\\
%\> is a $p\in \set{a,m,b}$ such that the three inclusion discs\\
%\> are contained in $D\paren{p, \frac{\rho_{k_p+1}(p)}{3}}$ then:\\
%  $k_a=k_m=k_b \le n$, and $\exists$ a $p\in \set{a,m,b}$ such that\\
% \> $\forall q \in \set{a,m,b}\sm \set{p}$, 
% $D(q, 3\rho_{k_q}(q)) \ib D(p, \frac{\rho_{k_p+1}(p)}{3})$ then:\\
3.a.\>\> $z_0 \as m$, $k \as k_m$, $g \as \fk$, $i \as 0$.\\
%\>\>\Comment{$p$ can be the point with the smallest $\beta_k(p)$ value}\\
4.\>\> While $\rho_k(z_{i}) \le 2^{5-2^{i}}\rho_k(z_0)$\\ % and $D(z_{i+1}, 3\beta_k(z_{i+1})) \cap I \neq \es$.\\
\>\>\>  $z_{i+1} \as z_i- g(z_i)/g'(z_i)$; $i \as i+1$.\\ %\Comment{Should we use the Schr\"oder iteration?}\\
\>\> $J \as  [z_{i-1} \pm 3 \rho_k(z_{i-1})]$\\%, z_{i-1} +3 \rho_k(z_{i-1})]$. \\
%\>\> \Comment{$J$ is the localization of the cluster; $i$ can be zero.}\\
5.\>\>If $w(J) \ge w(I)/2$ then return failure\\ %\Comment{The localization given by $J$ is worse compared to subdividing $I$}\\
%\>\>\> Return failure.\\
6.\>\>else return $(J, k)$.\\
7.\> Return failure.
}

We first explain some steps in the predicate above:
\begin{tightenum}{}
\item {\bf Step 2.} A point $p$ in $I$ can have more than one admissible value associated with it.
   The right admissible value is governed by $w(I)$,
  since we should only consider those clusters $\calc$ for which $r_\calc \ll w(I) \ll R_\calc$.
\item  {\bf Step 3.}  As $D(m, \rho_{k_m+1}(m)/3)$ contains all the three inclusion discs,
   they all contain the same cluster $\calc$. Otherwise, it is possible that the three
   inclusion discs contain different clusters but of the same size.
%   The condition states that $a$, $m$ and $b$ have the same admissible
%   value, say $k$; moreover, there is a point $p$ amongst them such all the three inclusion discs are in the 
%   disc $D(p, \rho_{k+1}(p)/3)$. Since this disc contains exactly $k$ roots, it
%   follows that all the three inclusion discs contain the same cluster $\calc$.
%   This rather strong condition will be satisfied in a neighborhood of 
%   a strongly separated cluster.\\
 \item {\bf Step 4.} This ensures that as $z_i$ converges to the root of $\fk$, the distance to $\calc$
  decreases quadratically; this fails when we are near $\calc$,
  or the root of $\fk$ is not near $\calc$. 
\item {\bf Step 5.} % It may so happen that $z$ is near a cluster $\calc$ to begin with (this will be the case, for instance,
  % after a successful call of the \texttt{Newton-Incl-Exc} predicate),  and hence the $\beta_k$
  % values do not decrease at all; in this case the localization  of $\calc$ given by $J$ is
  % not very useful, particularly when $w(J) \ge w(I)/2$, i.e., when subdividing $I$ gives a 
  % better localization than $J$.
  Required to ensure linear convergence to $\calc$. %Hence we declare failure, and go to subdividing $I$.\\
  % By now we know that $w(J) < w(I)/2$, and if $J \ibp I$ then we have a better localization of the cluster,
  % so we return $J$. \\
\item {\bf Step 6.}  The interval $J$ contains the cluster $\calc$. Moreover, 
  as $I \ib D(m, \rho_{k+1}(m)/3)$, we know that if the roots in $I$ are a
  subset of $\calc$, and hence are inside $J$. By now $w(J)< w(I)/2$,
  therefore, it suffices to return $J$. 
  % We do not yet know the location of $J$ with respect to $I$, 
  % which may cause  some issues (e.g., perhaps $\calc$ has been found earlier, 
  % or $J$ is not in $I_0$),  but these are handled by the parent routine.
\end{tightenum}
We now comment on some steps in \texttt{Newton-Isol}:
\begin{tightenum}{}
\item {\bf Step 1.a.} Ensures that a successful call to \texttt{Newton-Incl-Exc} is followed by a subdivision step.
   Thus the recursion tree is a binary tree. The predicate can still be successful on
   an interval $J$ returned by an earlier successful call.
  But the convergence in this case would only be linear, and so we prefer subdivision, though in practice
  one can go ahead with the linear convergence.
\item {\bf Step 2.b.} % We know that the interval $J$ returned by \texttt{Newton-Incl-Exc} contains a cluster $\calc$
  % and the roots in $I$ are from this cluster.
  Checks if  $\calc$ has not been found before (see \refLem{correctness}(iii)),
  and that $J$ is inside $I_0$;  if either of this test fails, then $I$ contains no roots and can be excluded.
\item {\bf Step 2.c.} As the only cluster in   $D(m(J), \frac{\rho_{k+1}(m(J))}{3})$
 is $\calc$, we can remove this disc from the intervals in $Q$.  It is this exclusion step that significantly
  contributes to the improvement of the subdivision algorithm.
\item {\bf Step 2.d.} This step adds the interval $J\si I_0$ containing the newly discovered cluster $\calc$
 to the set $\Phi$.% , but only the portion of inside $I_0$.
\end{tightenum}

There are only two loops in the algorithm: 
first, the while-loop in step 2 of the algorithm, and second, the Newton iteration
in step 4 of \texttt{Newton-Incl-Exc}. %We claim that both these loops will terminate.
The argument for the termination of the first loop is the same as for \texttt{Isolate}.
%  does not terminate, then there must be a sequence of intervals
% converging to a point $z$ such that neither $C_0$ nor $C_1$ hold for all these intervals; in the limit, this 
% means that a point is both a root of $f$ and not a root, which is a contradiction.
The termination of the second loop is guaranteed, because % for all iterates $z_i$ we know that the
% disc $D(z_i, 3\beta_k(z_i))$ contains a cluster of $k$ roots; therefore, $6\beta_k(z_i) \ge 2r_\calc$;
% moreover, as $f$ is square-free, the cluster radius $r_\calc > 0$; thus $\beta_k(z_i)$
% is bounded away from zero and hence the Newton-iteration has to terminate.
% Alternatively,
if $z_i$'s are such that
$\rho_k(z_i)$ keeps on decreasing, then in the limit $\rho_k$'s converge to zero; 
but the disc $D(z_i, 3\rho_k(z_i))$ contains exactly $k$ roots; 
since, in the limit $z_i$'s tend to a root $z^*$ of $f^{(k-1)}$, this implies that $z^*$ is a $k$-fold root 
of $f$, which is a contradiction as $f$ is square-free.

The following is a proof of correctness of the algorithm. 
\bthm
Given a polynomial $f$ and an interval $I_0$, \texttt{Newton-Isol}$(f, I_0)$ outputs a 
root partition $\calP$ of $I_0$.
\ethm
\bpf
We need to show the following claims:\\
1. $I_0\sm \calP$ contains no real roots of $f$.\\
2. $\calP$ contains (interior) pairwise disjoint intervals.\\
3. For all $I \in \calP$, $C_0$ or $C_1$ holds (follows from step 2.a.).

\refLem{correctness} gives us 
the correctness of \texttt{Newton-Incl-Exc}$(I)$, i.e., if the test is successful then it returns
 an interval $J$ such that any roots in $I$ are contained in $J$.
%  moreover, the only
% roots in the disc $D(m(J), \frac{1}{3\gamma_k(m(J))})$ are contained in the smaller
% disc $D(m(J), w(J)/2)$. These claims
% follows from \refLem{correctness} and the construction of
% the interval $J$ (note that $w(J)=3\rho_k(m(J))$). 
We only argue for the first claim.
For every interval $J$ returned by a successful call of the predicate, define
\beql{aj}
A_J \as D\paren{m(J) , {\rho_{k+1}(m(J))}/{3}} \sm D(m(J), w(J)/2),
\eeql
i.e., the annulus around $J$ that does not contain any roots. We exclude 
intervals if step (2.b) fails for the interval $J$, or a portion of an interval is removed in
step (2.c.). 
In the former case, either the cluster contained in $J$ was already detected, or it is outside $I_0$.
In the latter case, we do not loose any roots since $A_J$ has no roots.
So  $I_0\sm \calP$ contains no roots.
% The proof is by induction on the width of the input interval.
% The induction hypothesis is that the
% algorithm returns the correct output on intervals of width smaller than $w(I_0)/2$.
% The base case is when we do not have any recursive calls, i.e., 
% the \texttt{Newton-Incl-Exc} test always fails,  which implies that only $C_0$ and $C_1$ predicates are applied to the interval.
% But this is exactly the \texttt{Isolate} algorithm.
% % The set $\calI_0\as I_0 \sm \su_J \calA_J$, where the union is over all intervals
% % $J$ returned by a successful call to \texttt{Newton-Incl-Exc}. 
% Given the induction hypothesis, the three claims can now be
% shown to follow from our earlier observations about the interval $J$ 
% the following arguments:\\
% -- The annuli $A_J$ clearly do not contain any roots. Therefore, the real roots of $f$ in $I_0$
%   are  in fact in the set $I_0 \sm \su_J A_J$.\\
% -- The interior of the intervals in $P$ are disjoint; holds for the recursive
%   calls by the induction hypothesis. The intervals in $P$ are either obtained from the subdivided
%   process or by the recursive calls. Since the recursive calls are made on disjoint intervals 
%   (step 2.b ensures that), the interiors of the intervals obtained by recursion are disjoint.
%   The interior of intervals in $P$ obtained by subdivision process are trivially disjoint.\\
% -- For all intervals $I \in P$ either $C_0$ or $C_1$ holds. This is ensured by
%    step 2.a. of \texttt{Newton-Isol}.
\epf

% We also need the following observation:
% \bleml{obsv}
% A successful call to \texttt{Newton-Incl-Exc} predicate is followed by subdivision.
% \eleml
% \bpf
% Let $J=[z_i \pm 3\rho_k(z_i)]$ be the interval returned after a successful call.
% Since the while-loop in Step 4 stopped at $z_i$, it follows
% that 
% %$\rho_k(z_{i}) \le 2^{5-2^i}\rho_k(z_0)$ and 
% %$\rho_k(z_{i+1}) > 2^{5-2^{i+1}}\rho_k(z_0)$, which implies that 
% $\rho_k(z_{i+1}) > \rho_k(z_i)2^{-2^i} \ge \rho_k(z_i)/2$.
% If we call the predicate on $J$, the iteration will start from $z_i$
% and 
% \epf
\section{Complexity Analysis}\label{sec:complexity}
The main result is that \texttt{Newton-Incl-Exc} will be successful near a ssc $\calc$.
Let $c_0> 20$ be the constant in \refLem{correctness}, 
and $\calc$ a ssc throughout this section.
Our first claim is that $|\calc|$ is an admissible value for all points in $\calic$.
\bleml{con}
% Given a ssc $\calc$,
If $|z-m_\calc| \leq R_\calc/(8c_0n^2)$ then $\Delta_{k}(z) \ge c_0$.
\eleml
\bpf
Let $\alpha_1 \dd \alpha_k \in \calc$ and $\alpha_{k+1} \dd \alpha_n \in Z(f)\sm \calc$.
Moreover, assume that they are ordered in increasing distance from $z$.
From \refeq{ost}, we know that $2k|z-\alpha_{k+1}| > \rho_{k+1}(z) > |z-\alpha_{k+1}|/(2(n-k+1))$.
Moreover, $|z-\alpha_{k+1}| > R_\calc - |z-m_\calc| \ge R_\calc/2$; similarly, 
$|z-\alpha_{k+1}|< 3R_\calc/2$. Therefore,
        \beql{rk1}
        \frac{R_\calc}{4n} \le \rho_{k+1}(z) \le 3|\calc|R_\calc.
        \eeql
From \refeq{ost}, we again have $\rho_k(z) < 2k|z-\alpha_k|$. 
But $|z-\alpha_k| \le |z-m_\calc|+r_\calc$, which gives us
        \beql{rk2}
        \rho_k(z) \le 2k (|z-m_\calc|+ r_\calc).
        \eeql
Since $|z-m_\calc|, r_\calc \le R_\calc/(8c_0n^2)$, we get $\rho_k(z)\le kR_\calc/(2c_0n^2)$. 
Combining this with \refeq{rk1}, and the observation that $(n-k)k \le n^2/4$,
we obtain that $\Delta_k \ge 2c_0n^2/(8(n-k)k)\ge c_0$.
\epf

% We also need the following result on the location of the roots of the derivatives
% \bleml{derivrootsb}
% Given a strongly separated cluster $\calc$ of size $k$, for $j \le k$, there are 
% $k-j$ roots of the derivative
% $f^{(j)}(z)$ in $D(m_\calc, r_\calc)$ and the remaining $(n-k)$ roots are outside 
% $D(m_\calc, R_\calc/(2n^2))$.
% \eleml
Recall the definition of the intervals $I_\calc$, $\calic$ and the annulus $\annc$ from 
\refSec{notation}, and $A_J$ from \refeq{aj}.
\bleml{converse}
% Let $\calc$ be a ssc.
If an interval $I$ is such that 
        $$I \ib \calic=[m_\calc\pm R_\calc/(8c_0n^2)] \text{ and } w(I) > 72|\calc|r_\calc$$
then the pair $(J, k)$ returned by \texttt{Newton-Incl-Exc}$(I)$ is such that $k=|\calc|$,
$J \ib I_\calc=[m_\calc\pm 20kr_\calc]$, and $A_J\ip\annc$.
% returns the empty interval.
\eleml
\bpf
We show that the conditions on $I$ above imply that 
\texttt{Newton-Incl-Exc}$(I)$ reaches step 6 of \texttt{Newton-Incl-Exc} (all the steps below refer to the steps in the predicate).
This requires showing the following:
(i) all the conditions in step 3 are met; (ii) Newton-iteration in step 4
converges quadratically terminating with an interval $J$ with $w(J) < w(I)/2$,
and (iii) $J \ib I_\calc$. The following claims provide the proof.
% This will be done by
% showing that the admissible values computed in step 2 are well defined for the endpoints and
% the midpoint, and that conditions in step 3 hold for these values; this will be followed by showing that
Let $I=[a,b]$ and $m= m(I)$.

% We have already shown that for all $z \in I$, $\alpha_k(z) \ll 1$,
\begin{tightenum}{}
\item {\bf Claim 1:} For all $p \in \set{a, m , b}$, $k_p=|\calc|$. 
  Recall from Step 2 that $k_p$
  is defined as the {\em smallest} admissible value $k$ for which $I \ibp D(p, \rho_{k+1}(p)/3)$.
  From \refLem{con}, we have $k_p \le |\calc|$. 
  % Our claim is that the conditions of the lemma imply that this must be the case. 
  Since the roots in 
  $I$ can only come from $\calc$, any smaller admissible value corresponds to a subcluster $\calc'$ of $\calc$,
  which implies $R_{\calc'} \le r_\calc$. From \refeq{rk1} we know that
  $\rho_{|\calc'|+1}(p) \le 3(|\calc'|+1) R_{\calc'} \le 3|\calc|r_\calc$.
  Since $w(I) \ge 72|\calc|r_\calc$,  clearly $I\ibn D(p, \rho_{|\calc'|+1}(p)/3)$
  for any subcluster $\calc' \ibp \calc$. Thus $k_p \ge |\calc|$. \\

\item {\bf Claim 2:} For all $p \in I$,
   $I\ib D(p, \rho_{k+1}(p)/3)$. This will follow from the more general claim that
   $$D_1 \as D(m_\calc, R_\calc/(8c_0n^2)) \ib D(z, \rho_{|\calc|+1}(z)/3) \sa D_2,$$ 
   for all $z \in D_1$;   since $a, m, n \in I \ib D_1$, the claim holds.
   But for any $z \in D_1$, we know from \refeq{rk1} that
        $\frac{\rho_{|\calc|+1}(z)}{3} \ge \frac{R_\calc}{12n}$
 which is greater than $\frac{R_\calc}{4c_0n^2}$, the diameter of $D_1$, for $c_0 \ge 3$. 
%   The last quantity is the diameter of $D_1$.
%  therefore, $D_1 \ib D_2$. \\

\item {\bf Claim 3:} For all $z, w \in D_1$, the inclusion disc
   $D(z, 3\rho_k(z)) \ib D(w, \frac{\rho_{k+1}(w)}{3})$. This follows if
% This follows if we show that
%   $|z\pm 3\rho_k(z)-w| \le \frac{\rho_{k+1}(w)}{3}$,a
  % or if the following stronger claim holds:
   \beql{zw} |z-w| + 3\rho_k(z) \le \frac{\rho_{k+1}(w)}{3}. \eeql %\frac{1}{3\gamma_k(w)}.\eeql
   But $|z-w|, r_\calc \le R_\calc/(8c_0n^2)$, which along with \refeq{rk2} implies that
   $3\rho_k(z) \le 6kR_\calc/(4c_0n^2)$. Therefore, LHS of \refeq{zw} is smaller than 
   $13kR_\calc/(8c_0n^2)$, which is smaller than $R_\calc/(12n)$ for $c_0 \ge 20$, but
   from \refeq{rk1} we know that  the latter is smaller than the RHS of \refeq{zw}.\\

 \item {\bf Claim 4:} Let $z_i$ be the sequence of iterates computed in 
   the while-loop in Step 4. If $z_i \in D(m_\calc, \frac{R_\calc}{8c_0n^2})\sm D(m_\calc, 2r_\calc)$, 
   then $\rho_k(z_i)< 2^{5-2^i} \rho_k(z_0)$.  Since $z_i \nin D(m_\calc, 2r_\calc)$, $r_\calc \le |z_i - m_\calc|$,
   and hence from \refeq{rk2} we obtain  
   $\rho_k(z_i) \le 4k |z_i - m_\calc|$. From \cite[Thm.~2.2]{pawlowski:zeros-of-derivatives:99}
   we know that there is a unique root $z^*$ of $\fk$ in $ D(m_\calc, r_\calc)$.
   Therefore, $|z_i -m_\calc| \le |z_i - z^*|+r_\calc$. But as $z_i \nin D(m_\calc, 2r_\calc)$
   and $z^* \in D(m_\calc, r_\calc)$,
   we have $r_\calc \le |z_i- z^*|$, and hence $|z_i - m_\calc| \le 2|z_i - z^*|$.
   Thus, $\rho_k(z_i) \le 8k |z_i - z^*|$. As $z_0$ is an approximate zero to $z^*$ (see \refLem{correctness}(i)), 
   we know
   $|z_i - z^*| \le 2^{1-2^i} |z_0 - z^*|$, which implies that $\rho_k(z_i) \le 2^{4-2^i}k|z_0-z^*|$.
   Furthermore, from \refLem{correctness}(i) we know 
   $k|z_0-z^*|< 2\rho_k(z_0)$. Hence   $\rho_k(z_i)< 2^{5 - 2^i}\rho_k(z_0)$.\\

 \item {\bf Claim 5:} The interval $J \ib I_\calc$ and $w(J)< w(I)/2$. 
   The previous claim shows that if
   $z_i \nin D(m_\calc, 2r_\calc)$, then we will obtain quadratically decreasing values of $\rho_k(z_i)$.
   Thus when the iteration stops  $z_i \in D(m_\calc, 2r_\calc)$, and it follows from
   \refeq{rk2} that $\rho_k(z_i) \le 6kr_\calc$. Hence the interval
   $J = z_i \pm 3\rho_k(z_i)$ is contained in $I_\calc$, for $k \ge 2$.
   % If $I \ib \calA(m_\calc; kr_\calc, R_\calc/n^4)$ it follows that $I\cap J=\es$ and
   % hence we return an empty interval in Step 6.
   % If $I \cap D(m_\calc, kr_\calc) \neq \es$ and $w(I) \gg kr_\calc$, then from the
   Moreover, $w(J) \le 36kr_\calc < w(I)/2$, and hence the condition in 
   Step 5 fails and we return $J$. The claim on the annulus follows from \refeq{rk1}.
\end{tightenum}
\epf

The following result translates the result above in terms of the subdivision tree: 
\bcorl{con}
Let $\calc$ be a ssc such that $\calic \ib I_0$.
If $I$ is the first interval such that 
\texttt{Newton-Incl-Exc}$(I)$ is successful and the interval returned contains $\calc$, then
$\calic \ib I' \su I''$, where $I'$ is the parent-interval of $I$
and $I''$ is one of $I'$'s neighbors.
\ecorl
\bpf
In the worst case, $\calc$ will be detected 
the first time in the subdivision tree an interval $I\ib\calic$.
% then \texttt{Newton-Incl-Exc}$(I)$ is successful. 
For such an $I$, we show $w(I) \gg kr_\calc$. 
Since $I$ is the first interval to fall in $\calic$, both
$I'$ and $I''$ have endpoints outside $\calic$, thus $\calic \ib I' \su I''$.
% the width of the interval associated with parent $I'$ of $I$
% and one of the neighbors of $I'$ cover $\calic$ completely.
So $2w(I) \ge R_\calc/(16c_0n^2)>72 kr_\calc$, as $\calc$ is ssc.
%\footnote{The constant involved in ssc is $16c_0\times 72$.}
The claim clearly holds if  $\calc$ is detected at an ancestor of $I$. 
%  one of the endpoints of $I'$ and $I''$ is outside $\calic$.
% If either one of them, say $I'$, is contained in $\calic$, then \texttt{Newton-Incl-Exc}$(I')$ is successful.
% Let $I'$ be the interval associated with
% the parent node of $I$. Since the condition in step 5 failed for $I'$ wrt $\calc$,
% we know that one of the endpoints of $I'$ is outside the disc $\calic$.
% There are two cases to consider:
% \begin{enumerate}
% \item If the right endpoint of $I'$ is outside $\calic$. In that case, $I$ is the left child of $I'$.
%   Moreover, the interval $I''$ to the left of $I$, which has the same width as $I$,
%   would have failed the \texttt{Newton-Incl-Exc} test, since otherwise
%   $I$ would have been removed from the interval. Thus the left endpoint of $I''$ is outside $\calic$.
%   Therefore, $3w(I) \ge R_\calc/n^4$, and since $\calc$ is a strongly separated cluster it follows
%   that $w(I) \gg kr_\calc$.
% \item If the left endpoint of $I'$ is outside $\calic$, then $I$ is the right child of $I'$. 
%   Let $I''$ be the interval to the right of $I'$ in the subdivision tree. Again the \texttt{Newton-Incl-Exc} test must
%   have failed for $I''$, and hence its right endpoint is outside $\calic$. Therefore,
%   $4w(I) \ge R_\calc/n^4 \gg kr_\calc$.
% \end{enumerate}
% Note the asymmetry of the two cases; this arises because we are doing a bfs of the subdivision tree.
\epf
%%%%%%%%%%%%%%%%%%%%%%%%%%%%%%%%%%%%%%%%%%%%%%%%%%
%\subsection{Bounding the Subdivision Tree}

\Remark The proof above gives us the explicit constant in the definition of ssc, namely, we require
$R_\calc/r_\calc > 16c_0\times 72 n^3$. A careful working out of the proofs shows that 
the weaker inequality $R_\calc/r_\calc > 4c_0 \times 72 (n-|\calc|)|\calc|^2$, (or even
$50 c_0n^3$) is sufficient.
% or even $R_\calc/r_\calc > 42 c_0 n^3$.

Recall that the set of all roots $Z(f)$ is a cluster.
As a consequence of \refLem{converse}, we assume that $I_0\ib n I_{Z(f)}$; 
otherwise $\texttt{Newton-Incl-Exc}$ will be successful right away and the 
interval returned will satisfy the property.

\subsection{An Integral Bound on the Subdivision Tree}
Let $\calno$ be the set of leaves in the subdivision tree of \texttt{Newton-Isol}$(f,I_0)$.
Step 1.a. of the algorithm ensures that
the subdivision tree is a binary tree. Therefore, it suffices to bound $|\calno|$.
For this purpose, we use the general framework of continuous amortization developed in
\cite{burr-krahmer-yap:continuousAmort:09} and generalized in \cite{burr:contamort:13}.
The idea is to bound $|\calno|$ by an integral and then derive an upper bound on this integral.
For this purpose we need the following notion: Given a choice of predicates $C_0$, $C_1$, a 
function $G:\RR\to\RR_{\ge 0}$ is called a \dt{stopping function} corresponding to $C_0$ and $C_1$ 
if for every interval $I$, if there is an $x\in I$ such that 
$w(I)G(x) \le 1$, then either $C_0(I)$ or $C_1(I)$ holds.
Stopping functions, corresponding to different predicates, are 
provided in \cite{burr:contamort:13}. The crucial property of $G(x)$ is the following:
\bleml{cp}
If $C_0(I)$  and $C_1(I)$ fail for an interval $I$, then for all $J \ib I$, such that
$2w(J) \ge w(I)$, $2\int_{J} G(x) dx \ge 1$. 
\eleml
\bpf
From the definition of $G(x)$, we have
for all $x\in I$, $G(x) w(I) \ge 1$. As $J \ib I$, $\forall x\in J$, $2G(x)w(J) \ge G(x)w(I) \ge 1$.
Thus $2\int_{J} G(x) dx \ge 2 w(J) \min_{x \in J}G(x) \ge 1$.
\epf

% If $C_0(I')$  and $C_1(I')$ both fail, then for all $x\in I',$ $G(x) w(I') \ge 1$, 
% and hence $\int_{I'} G(x) dx \ge 1$. Moreover,

% If $I$ is a leaf and $I'$ its parent, then we have for all $x \in I$, $G(x)w(I')=2G(x)w(I) \ge 1,$ and hence
% at the leaves, $2 \int_IG(x) dx \ge 1$ (note the crucial property that $w(I')=2w(I)$).

The main result of this section is the following:
\bthml{intbound}
        $$|\calno| \le 4n + 2\int_{I_0 \sm \su_\calc \annc} G(x) dx,$$
where the union is over all ssc $\calc$ in $T_f$.
\ethml
We bound $\calno$ recursively.
The leaves in $\calno$ correspond to three types of intervals:
\begin{tightenum}{r}
\item  intervals in the root partition $\calP$,
\item  intervals that were discarded in step 2.c., and
\item  intervals for which condition 2.b fails to hold (either cluster already found, or $J\si I_0=\es$).
\end{tightenum}
We will bound each of these three types. 
We analyse what happens before the first set of recursive calls.

Let $\Phi$ be the set of intervals collected in Step 2.d.~of the algorithm,
$A_J$ be as defined in \refeq{aj}, and $\cali_J \as J\su A_J$. 
From the construction of $\Phi$, we know that all intervals $J \in \Phi$ are contained in $I_0$
and each contains a unique cluster.
For each $J \in \Phi$, let $L_J$ be the set of {\em parent-intervals} of intervals in $Q$ 
that intersect $\cali_J$; the type (ii) intervals are children of intervals in $L_J$. Let
$M_J$ be the set of intervals that do not intersect $\cali_J$ and are of type (iii).
See \refFig{3types} for an illustration of these types.
% A type (i) or type  (iii) interval could be a subset of an interval in $I \in L_J$;
Note that if $I\in L_J$ contains an endpoint of $\cali_J$, then
$I\sm \cali_J$ can be of type (i) or (iii); but there can be
at most two such intervals for each $J$ on either side of $\cali_J$.
We abuse notation and use $L_J$ to represent a set as well as the union
of the intervals in it; same for $M_J$.

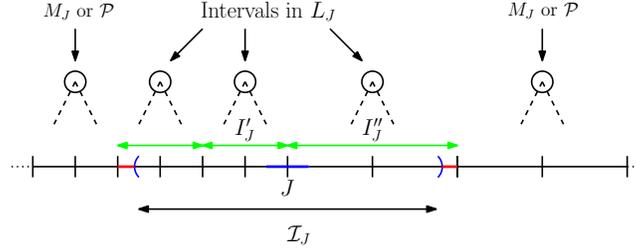
\vfigpdf{The three types of intervals in $\caln(I_0)$. Intervals in $L_J$ are shown in green.
  The remaining intervals could be in $M_J$ or $\calP$.
  The width of the red colored intervals  can be much smaller than their parents.
  But there are at most two such intervals.}{3types}{0.5}

% We first bound $|M_J|$.
For an $I\in M_J$, both $C_0$ and $C_1$
failed. Therefore, from \refLem{cp} we get
$|M_J| \le 2\sum_{I \in M_J}\int_I \gxdx= 2\int_{M_J}\gxdx$.
As the predicates $C_0$ and $C_1$ also fail for the intervals in $L_J$,
we can similarly bound $|L_J|$. But this 
effectively amounts to doing subdivision on $J$.
% annuls the advantage of doing recursion on $J$.
To avoid this we do the following: since the width of the intervals in $L_J$ is 
more than $w(J)$, we know that there are at most two neighboring
intervals $I'_J$ and $I''_J$ that contain $J$. We count them separately, and for the rest
we use \refLem{cp} to get
$|L_J|\le 2 + 2\int_{L_J \sm (I'_J \su I''_J)} G(x) dx$. 
For an interval  $I \in \calP$, we expect $2\int_IG(x) dx \ge 1$, as the predicates
must have failed for the parent $I'$ of $I$. However, 
\refLem{cp} requires that $w(I')\le 2w(I)$. 
This can fail to happen near the boundary of $\cali_J$, as noted earlier. 
But then there are at most two such intervals.
Therefore, the number of intervals in $\calP$ coming from the non-recursive calls is at most
$2|\Phi| + 2\int_{I_0 \sm \su_J (L_J \su M_J)} \gxdx$.
Combining this with the bounds on $|L_J|$ and $|M_J|$ we get
\beql{ni1}
|\calno| \le 4|\Phi| + 2\int_{I_0 \sm \su_J(I'_J \su I''_J)} \gxdx + \sum_{J \in \Phi}|\caln(J)|.
\eeql
To open the RHS recursively, we introduce the notion of \dt{cluster tree $\calto$
with respect to an interval $I_0$}: It is the smallest subtree $\tcc$ of $T_f$ rooted at
a cluster $\calc$ such that $I_0\ib I_\calc$; since by assumption
$I_0 \ib nI_{Z(f)}$, in the worst case, $\calto$ is $T_f$.
Moreover, as enlarging $I_0$ increases the integral in \refeq{ni1}, we further
make the simplifying assumption that $I_0=2I_{\calco}$,
where  $\calco$ is the root of $\calto$. 

Let $\calc$ be the cluster associated with a node $u$ in $\calto$. 
Let $J_u \in \Phi$ be the interval returned the first time 
$\calc$ is detected by \nie. 
Define $A_u \as (I'_{J_u} \su I''_{J_u}) \sm J_u$; if $\calc$ is not detected,
let $A_u = J_u = \es$.
% With a node $u$ in $\calto$ we associate an interval $J_u\as J$ 
% if the interval $J$ containing the cluster associated with $u$ (note that $J$
% may or may not intersect $I_0$) causes a recursion of the algorithm; 
% otherwise $J_u \as \es$; if $J_u \neq \es$ then further 
% define $A_{u} \as (I'_{J_u} \su I''_{J_u}) \sm J_u$.
Using this notation, the following bound can be derived from \refeq{ni1} by induction:% opening the RHS:
%\bleml{bound}
        \beql{ni2}
        |\calno| \le 4|\calto| + 2 \int_{I_0 \sm \bsu_{u \in \calto} A_u}G(x) dx.
        \eeql
%\eleml
% This bound is not very useful as it 
% depends on the intervals $J$ computed by the algorithm. 
For a ssc $\calc\in \calto$, the assumption $I_0=2I_{\calco}$ ensures that
$\calic\ib I_0$. So \refCor{con} implies that $I'_u \su I''_u\ip \calic$, and
\refLem{converse} implies that $J_u\ib I_\calc$; hence, $A_u \ip \annc$.
Considering only the ssc in $\calto$ on the RHS of \refeq{ni2} we obtain 
% We want something that can be
% bounded {\it a priori} in terms of $I_0$. Consider a node $u$ corresponding to a ssc $\calc$
% in $\calto$. From \refCor{con} it follows that $A_u \ip \annc$.
% Therefore, removing $\annc$ from the RHS of \refeq{ni2} we obtain 
% From \refCor{con}, we know that
% for any ssc $\calc$ in $\calto$ such that $\calic$ intersects $I_0$
% will be found by the \texttt{Newton-Incl-Exc} predicate; the exceptions are those ssc $\calc$ for which
% $w(I_0 \si \calic) \ll kr_\calc$, but there can be at most two such clusters;
% the contribution of such clusters will be bounded later.
% Moreover, the annuli rejected for such a $\calc$ contains $\annc$.
% Therefore, by keeping only the terms in the integral above corresponding to
% the ssc in $\calto$, the value of integral only increases.
% For simplicity, however, we can consider all strongly separated clusters in $T_f$,
% to obtain the following:
% \bleml{bound}
        \beql{nie3}
        |\calno| \le %2|\calto| + 2 \int_{I_0 \sm \su_{\calc \in \calto} \annc}G(x) dx
                          4|\calto| + 2 \int_{I_0 \sm \su_{\calc} \annc}G(x) dx,
        \eeql
% where $\calc$ is any ssc in $\calto$.
As $|\calto|\le n$, we get \refThm{intbound}.

\subsection{Bound for the Descartes's rule of signs}\label{sec:ndsc}
% So far the argument is independent of the choice of the stopping function.
In this section, we derive the following bound:
\bthml{bound}
Given a square-free polynomial $f\in \RR[x]$ of degree $n$, 
the size of the subdivision tree constructed by \texttt{Newton-Isol}($f, I_0$) using
predicates based on the Descartes's rule of signs is bounded by $O(n\ln n)$.
\ethml
% Most of the analysis will be independent of the choice of the predicate.
% It is only towards the end that we actually use th

We bound the RHS of \refeq{nie3},
%$\int_{I_0 \sm \su_{\calc } \annc} G(x) dx$.
where the stopping function corresponds to the Descartes's rule of signs. 
We use the same stopping function as described in \cite{burr:contamort:13},
but explain why the argument there fails to give us the bound above.% in \refThm{bound}.

Let $V \as Z(f)$, the set of roots of $f$. Define  $d(x, V)$ as the distance from
$x$ to the closest point in $V$, and $d_2(x,V)$ as the distance to the second closest point in $V$. 
The crucial idea in \cite{burr:contamort:13} is to partition the integral
over the (real) Voronoi interval $I_\alpha$ of each root $\alpha$ (for the moment
suppose $\alpha \in \RR$).
% For a root $\alpha \in \RR$, let $\sigma_\alpha \as d_2(\alpha, V)$.
%and if $\alpha\in \CC\sm \RR$, then $\sigma_\alpha \as |\Im(\alpha)|$.
% The stopping function is $2/d_2(x,V)$ in a certain neighborhood of $\alpha$, and $1/d(x,V)$ 
% else where. More precisely, 
Define $J_\alpha\as [\alpha\pm \frac{d_2(\alpha, V)}{2}]$. Then
for $x\in J_\alpha$, $G(x) \as 2/d_2(x,V)$, and for $x\in I_\alpha\sm J_\alpha$,  $G(x) \as 1/|x-\alpha|$.
Break $\int_{I_\alpha}\gxdx$ as $\int_{J_\alpha}\gxdx + \int_{I_\alpha\sm J_\alpha}\gxdx$.
In \cite{burr:contamort:13} it is shown that the first integral is $O(1)$,
% $\int_{J_\alpha} \gxdx = O(1)$, 
and the second integral is $O(\log w(I_\alpha)/d_2(\alpha, V))$;
%  $\int_{I_\alpha\sm J_\alpha} \gxdx = O(\log w(I_\alpha)/d_2(\alpha, V))$;
from Cauchy's bound we can assume that $w(I_\alpha)=2^{O(L)}$.
The problem is that
in the worst case this ratio can be $\Omega(n(L+\log n))$. E.g., if all the other roots 
are of the form $\alpha \pm i t$, for increasing values of $t$, then $I_\alpha$ is the x-axis.
Therefore, 
$\int_{I_\alpha\sm J_\alpha}\gxdx= \Omega(L - \log d_2(\alpha, V))$; in the worst case $d_2(\alpha, V)$ can be 
the root separation bound.

Our idea is based on the observation that roots with very small separation give rise to root clusters.
For clusters that are not ssc, the ratio $R_\calc/r_\calc=O(n^3)$, therefore, the number of subdivisions
needed to bridge this gap is $O(\log n)$. For a ssc, the gap is bridged by \texttt{Newton-Incl-Exc}
so that the subdivision is restricted to the ranges $R_\calc$ to roughly $R_\calc/n^2$ and
$|\calc|r_\calc$ to $r_\calc$, both of which take $O(\log n)$ subdivisions.
Doing this for all clusters basically gives the bound in \refThm{bound}.
% So what matters in the argument above is the relative ratio of the width of the cluster to the
% minimum separation of the roots inside the cluster. For clusters that are not ssc, this ratio is roughly 
% bounded by $O(n^3)$. For ssc, \texttt{Newton-Incl-Exc} helps us detect them and 
% reduce this relative ratio to $n^{O(1)}$. % This us to the result in \refThm{bound}.

Let $P\ib \CC$ be a pointset such that any non-real point in $P$ also has its complex conjugate in $P$.
Such a set of points is called \dt{dense} if no proper subset of $P$ forms a non-trivial cluster, i.e.,
for all $S \ibp P$, such that $|S| \ge 2$,
the disc $3D_S$ contains a point from $P\sm S$. This structure plays a fundamental
role in our arguments, as do the following two integrals (see \cite{burr:contamort:13,sharma-yap:near-optimal:12}):

\bleml{int}
	Let $\gamma\in \CC$ and $J=[r,s]$. 
	\\{\bf(Re)} If $\gamma\in \RR\sm J$, then
		\beql{areal}
		\int_J \frac{dx}{|\gamma-x|} = \ln
                \left| \frac{\gamma-s}{\gamma-r}\right|^{\delta(J>\gamma)},
		\eeql
         where $\delta(J> \gamma)= +1$ if $r > \gamma$ and $-1$ if $s < \gamma$.
	\\{\bf(Im)} If $\gamma \in \CC\sm \RR$ then
        %, say $\gamma= \Re(\gamma) + \ii \Im(\gamma)$, 
        % define 
       %          \beql{rho}
       %          \rho_t \as \frac{2|t-\gamma|}{|\Im(\gamma)|},
       %          \eeql
       % where $t\in\set{r,s}$; note that $\rho_t \ge 1$. Then
		\beql{acomp}
%                \begin{split}
		\int_J \frac{dx}{|\gamma-x|} 
                                \le 
                                O\paren{\ln \frac{\max \set{|s-\gamma|,|r-\gamma|}}{|\Im(\gamma)|}}.
                                % \begin{cases}
                                % \ln 4 2 \frac{|s-\gamma||r-\gamma|}{|\Im(\gamma)|^2}&\text{if $\Re(\gamma) \in J$}\\
                                % \ln 2\paren{\frac{|\gamma-s|}{|\gamma-r|}}^{\delta(J>\Re(\gamma))} &\text{otherwise}.
                                % \end{cases}
%                                \frac{\rho_s^{\delta(s\ge\Re(\gamma))}}{\rho_r^{\delta(r>\Re(\gamma))}}.
%               \end{split}
		\eeql
          % Note that when $\Re(\gamma) \in [r, s]$ then both the terms are in
          % the numerator.
\eleml

We now give the proof of \refThm{bound}.\\
\bpf
The proof is by induction on $|\calto|$.
We claim that 
        \beql{claim1}
        \int_{I_0 \sm \su_{\calc} \annc } \gxdx = O(|\calto|\ln n).
        \eeql
Let $\calco$ be the root of $\calto$; by assumption we have $I_0= 2I_{\calco}$. 
% Since choosing a larger interval increases the value of the integral,
% we  can assume that the interval $I_0=|\calco|I_{\calco}$. 
Let $\calmo$ be the children of $\calco$ in $\calto$.
% A \dt{maximal ssc} $\calc$ in $\calto$ is a ssc for which there is {\em no} ssc $\calc'$
% such that $\calc \ibp \calc' \ibp \calco$.
% Let $\calM_0$ be the set of maximal sscs in $\calto$. 
% Then it is clear that for any other ssc $\calc'$ in $\calto$,
% the interval $I_{\calc'}\ib I_\calc$, for some $\calc \in \calmo$.
Consider a ssc $\calc \in \calmo$. Then $I_0 \sm \annc \ib (I_0 \sm \calic) \su 2I_\calc$.
If $\calc'$ is a ssc contained in $\calc$, we can inductively remove $\calA_{\calc'}$ from $I_\calc$.
This also works for clusters that are not ssc in $\calmo$, since by definition $\calic=2I_\calc$.
Therefore,
        $$I_0 \sm \su_{\calc} \annc 
                \ib (I_0 \sm \su_{\calc \in \calmo} \calic) \su \paren{\su_{\calc \in \calmo}\paren{2I_\calc \sm \su_{\calc'\ibp\calc} \calA_{\calc'}}}.$$
% That is, first remove from $I_0$ the intervals $\calic$ corresponding to the maximal sscs, and 
% recursively do the same for each maximal ssc. 
We claim that 
        \beql{claim2}
        \int_{I_0 \sm \su_{\calc \in \calmo} \calic}G(x) dx = O(|\calmo|\ln n).
        \eeql
As $|\tcc| < |\calto|$, for $\calc \in \calmo$, by induction we obtain
        $$\int_{2I_\calc \sm \su_{\calc'\ibp\calc} \calA_{\calc'}}\gxdx = O(|\tcc|\ln n).$$
This bound along with \refeq{claim2} and the
observation that $|\calmo|+\sum_{\calc \in \calmo}|\tcc|<|\calto|$ gives us \refeq{claim1}.
The base case is when $\calmo$ contains only leaves, 
in which case \refeq{claim1}  reduces to \refeq{claim2}.

We next claim that
        $$
        \int_{I_0 \sm \su_{\calc \in \calmo} \calic}\gxdx= O(\ln n) + 
                        \int_{I'_0 \sm \su_{\calc \in \calmo}\calic} \gxdx,
         $$
where $I'_0 \as [m_{\calco} \pm 2r_{\calco}]$.
If $\calco$ is not a ssc, then this is clear as $I'_0 =2I_{\calco}=I_0$. %and the RHS is only larger.
If $\calco$ is a ssc, then $I_0 = [m_{\calco}\pm |\calco| r_{\calco}]$. Break $I_0$ as 
$I'_0$, $[m_{\calco}+2r_{\calco}, m_{\calco}+ |\calco|r_{\calco}]$
and $[m_{\calco}-2r_{\calco}, m_{\calco}- |\calco|r_{\calco}]$. 
%We claim that the integral on the latter two intervals is bounded by $O(\ln |\calco|)$. 
The closest root to any $x$ in these intervals
is from $\calco$. Moreover, as $|x-m_{\calco}| \ge 2r_{\calco}$, we get
$G(x) \as 1/d(x, V) \le 2/|x-m_{\calco}|$. Therefore, from \refLem{int}(Re) it follows that
$\int_{m_{\calco}+2r_{\calco}}^{m_{\calco}+ |\calco|r_{\calco}} \frac{2}{|x-m_{\calco}|} = O(\ln |\calco|)$.
Similarly for the other interval. Hence to prove \refeq{claim2}, it suffices to show
        \beql{claim3}
        \int_{I'_0 \sm \su_{\calc \in \calmo}\calic} \gxdx = O(|\calmo|\ln n).
        \eeql

Let $\calmo$ also denote the pointset obtained
by replacing each $\calc \in \calmo$ by its center $m_\calc$. 
% , effectively
% making the children of $\calco$ all leaves.
We will use \refLem{dense} to prove \refeq{claim3}.
As no subset of $\calmo$ forms a cluster, $\calmo$ is a dense pointset, and \refLem{dense} is applicable.
However, we first remove some region around every $p \in \calmo\si \RR$ to be able
to invoke \refLem{dense}. For each such $p$, define
$J_p \as [p \pm d_2(p, \calmo)/2]$. If  $p=m_\calc$, for $\calc \in \calmo$, then
$\calip \as \calic \ib J_p$. We claim 
        \beql{claim4}
        \int_{\su_{p\in \calmo} (J_p \sm \calip)}\gxdx = O(|\calmo|\ln n).
        \eeql
From \refLem{dense} we get
        $$\int_{I'_0 \sm (\su_{p \in \calmo}J_p)}\gxdx = O(|\calmo|\ln n).$$
Combining these two bounds, along with the observation that the union of the sets
$I'_0 \sm (\su_{p \in \calmo}J_p)$ and $\su_{p\in \calmo} (J_p \sm \calip)$ is the set
$I'_0\sm \su_{p\in \calmo} \calip$, completes the proof of \refeq{claim3}.

To prove \refeq{claim4}, we show that $\int_{J_p\sm \calip}\gxdx= O(\ln n)$, and then sum over all
$p \in \calmo$. There are three cases to consider: 
\begin{tightenum}{r}
\item  $p=m_\calc$ for some normal cluster $\calc \in \calmo$.
%The integral over $2I_\calc$ has been done inductively.
Then $J_p= [m_\calc \pm R_\calc/2]$ and $\calip=\calic=[m_\calc\pm 2r_\calc]$.
Therefore, $J_p \sm \calip$ contains $[m_\calc+2r_\calc, m_\calc+R_\calc/2]$
and $[m_\calc-R_\calc/2, m_\calc - 2r_\calc]$. 
The nearest root to any $x$ in these two intervals is from $\calc$. 
Since $x$ is outside $2I_\calc$, it follows that $d(x,V) \ge |x-m_\calc|/2$.
Therefore, $G(x)\as\frac{1}{d(x,V)} \le 2/|x-m_\calc|$. 
From \refLem{int}(Re), we obtain
$\int_{m_\calc+2r_\calc}^{m_\calc+R_\calc/2}\gxdx=O(\ln R_\calc/r_\calc)$.
Since $\calc$ is not a ssc, $R_\calc/r_\calc = O(n^3)$, which gives us the desired bound.
The same applies to the other interval. % $[m_\calc-R_\calc/2, m_\calc - 2r_\calc]$. 
\item Suppose $p=m_\calc$, where $\calc \in \calmo$ is a ssc. 
%The integral over ${\calic}$ has been done inductively. 
Then $J_p = [m_\calc \pm R_\calc/2]$ and $\calip=\calic=[m_\calc\pm R_\calc/n^2]$.
% Consider the integral on $J_p \sm \calic$.
Let $I'_\calc \as [m_\calc +\frac{R_\calc}{n^2}, m_\calc + \frac{R_\calc}{2}]$ be one of the  
intervals in $J_p \sm \calip$. The nearest root to any $x\in I'_\calc$ is from $\calc$. 
Since $x\nin 2I_\calc$, it follows that $d(x,V) \ge |x-m_\calc|/2$.
Therefore, $G(x)\as\frac{1}{d(x,V)} \le 2/|x-m_\calc|$. 
Applying \refLem{int}(Re), we get $\int_{I'_\calc} \gxdx \le 4\ln n$.
Similarly, for the other interval.
%$[m_\calc -\frac{R_\calc}{2}, m_\calc - \frac{R_\calc}{n^2}]$. 
\item  $p$ is a real root then $\calip=\es$. 
For $x \in J_p$, our stopping function $G(x)=2/d_2(x,P)$, i.e., corresponding
to the inclusion predicate. Suppose $q \in P$ is such that $d_2(p, P)=|p-q|$.
Then for all $x \in J_p$, $d_2(x,P) \ge |p-q|-|p-x|\ge \sigma_p/2$,
and hence 
        $\int_{J_p} \frac{2dx}{d_2(x,P)} \le 4\int_{p-\sigma_p/2}^{p+\sigma_p/2} \frac{dx}{\sigma_p}=O(1)$.
\end{tightenum}
\epf

The proof above can carried out with the exact constants involved in the definitions of 
$I_\calc$, $\calic$ and $\annc$ (see \refLem{converse}), 
but they will be absorbed by the big-O notation. Note that the $O(n \ln n)$ bounds
the number of calls to the $C_0$ predicate.
The specialization of $G(x)$ for $C_0$ is  $1/d(x,V)$.
The corresponding specialization for Sturm sequences is $1/d(x, V\si \RR) \le 1/d(x,V)$. 
Therefore, $O(n\ln n)$ holds for \texttt{Newton-Isol} combined with Sturm sequences. For \texttt{Eval},
one specialization of the stopping function for the $C_0$ predicate is $n/d(x,V)$, 
which immediately gives an $O(n^2\ln n)$
bound for \texttt{Newton-Isol} combined with \texttt{Eval}. Whether it can be improved
using the more precise specialization $\sum_{\alpha\in V}\frac{1}{|x-\alpha|}$ remains open.

Let $P$ be a dense pointset $P$.
Given a point $p \in P$, define $\sigma_p \as \min|p-q|$, where $q\in P\sm \set{p}$, 
and $J_p \as [p\pm \sigma_p/2]$.
We want to bound $\int_{(2D_P \si \RR) \sm \su_p J_p }dx/d(x,P) $.
We first show an $O(|P|^2)$ bound, essentially following \cite{burr:contamort:13}.
Let $\calv_p$ be the set of points in $2D_P \si \RR$ closer to $p$ than to any other point in $P$. 
% , i.e., the (real) Voronoi region of $p$. 
It is clear that $J_p \ib \calv_p$.
The intervals $\calv_p$ partition $2D_P \si \RR $.
Then $\int_{\calv_p\sm J_p} dx/d(x,P)$ can be shown to be bounded by
 $O(\ln (r(D_p)/\sigma_p))$. Using the density of $P$, it can be shown
that if $p, q$ are such that $\sigma_p = |p-q|$ then $P \ib 3^{O(|P|)}D_{\set{p,q}}$,
which implies that  $r(D_P) \le 3^{O(|P|)}\sigma_p$,
for all $p \in P$. 
% For any point $p$, let $q \in P\sm \set{p}$ be a closest
% point to it. Then from the density of $P$, it follows that $3D_{\set{p, q}}$ 
% contains another point $s$; now applying the same argument to $\set{p,q,s}$,
% we get a new point in $3^2D_{\set{p,q,s}}$;
% continuing in this manner we can cover all the points in $P$ in $3^{|P|}D_{\set{p,q}}$.
% Thus for all $p \in \callc$, $\int_{I_p}\gxdx = O(|P|)$ and hence
% $\int_{2D_P\si \RR}\gxdx = O(|P|^2)$. 
This gives an $O(|P|^2)$ bound instead of the bound in \refThm{bound}.
To obtain that we need to amortize the integral carefully.
The intuition is that if $\sigma_p$ is very small then there there must a lot of other points
close to $p$, and hence the width of $\calV_p$ cannot be very large compared to $\sigma_p$ .
% Moreover, we cannot have $|P|/2$ such pairs, as the $3^{|P|}$
% inflations of all the corresponding discs have to cover $P$, 
% and hence there must be considerable overlap amongst these discs.
The challenge is to get an ``almost cluster-like'' decomposition of $P$.
We construct a tree on $P$ that gives us this decomposition.

We describe an iterative bottom-up procedure to construct a tree $\caltp$ with leaves from $P$.
% For $p \in P$, recall that $\sigma_p$ is the distance to the closest point from $p$ to $P\sm \set{p}$.
% If $p \in \callc \si \RR$,  let $D_p \as D(p, \sigma_p/2)$;
% if $p \in \callc \sm \RR$, then let $D_p\as D(\Re(p), \Im(p))$.
% Order the $\sigma_p$'s in increasing order, and let  
% $\sigma_1 < \sigma_2  < \cdots < \sigma_\ell$ be the distinct values of these distances.
% At step $i$, we have a partition $G_1 \dd G_{m}$ of $\callc$ into components (initially they are 
% the individual points $p$). Each component $G_j$ has an associated region 
% $\rho_{G_j}$
Let $\sigma \as \min_{p \in P}\sigma_p$.
For all points $p \in P$, draw a disc of radius $\sigma/2$ centered at $p$.
As $\sigma$ is the smallest distance between any pair of points,
two such discs can  at most touch each other. The discs touching each other form a connected
component. The collection of the largest connected components partitions $P$ (leaves are 
considered as components).
Moreover, there is at least one component $G\ib P$ that has cardinality strictly greater than one;
the components with cardinality one are the leaves.
% For each component $G \ib \callc$, let $\rho_G\as \su_{p \in G}D_p$. 
% If the regions corresponding
% to  two components overlap, then we merge the components into a single component. 
For each such component $G$,
we introduce an internal node $u$ in $\caltp$ with children as the leaves $p$, where $p \in G$;
let $G_u \as G$, the associated component, and $\sigma_u \as \sigma$.
% $\rho_u$ the associated region and $\sigma_u\as \sigma_1$. It is easy to see that
% if $\Delta_u$ is defined as the least distance from $p \in G_u$ to a point in $\callc\sm G_u$
% then $\Delta_u > \sigma_u$.
% From construction it is clear that the regions corresponding
% to the components are disjoint and at distance at least $\sigma_1$.
Now redefine $\sigma$ as the minimum separation between the 
%  repeat the step with $\sigma_2 \dd \sigma_\ell$ with the 
components constructed so far,
draw a disc of radius $\sigma/2$ centered at each $p\in P$, and continue as above.
Let $\caltp$ be the tree constructed in this bottom-up manner; see \refFig{ripples}.
% , i.e., further 
% increase the radius of the disc around every point to $\sigma_2$.
% Again consider the largest connected components formed from the existing
% components and leaves; for every new component that comes into existence (there is at least one), introduce a
% corresponding node $u$ in $\caltp$, and assign $\sigma_u \as \sigma_2$. 
% We continue in this manner with $\sigma_3 \dd \sigma_\ell$.
% For every node $u \in \caltp$, 
% $\sigma_u$ captures the separation $\sigma_i$ that formed the component.
% In addition to $G_u$ and $\sigma_u$, 
Further define the following quantities for each $u \in \caltp$:
\begin{tightenum}{r}
\item  $\nu_u$ as the number of children of $u$,
\item $m_u$ be the center and $r_u$ be the radius of $D(G_u)$. %, and\\
\end{tightenum}
% (iii) $\Delta_u$ as the least distance $|p-q|$, where $p \in G_u$ and  $q\in P\sm G_u$.
% (iii) The region $\rho_u \ib \RR^2$ as the union of the discs
% $D(p, \sigma_u/2)$, when $p \in G_u \si \RR$, and $D(\Re(p), \Im(p)+\frac{\sigma_u}{2})$, if $p\in G_u \sm \RR$.

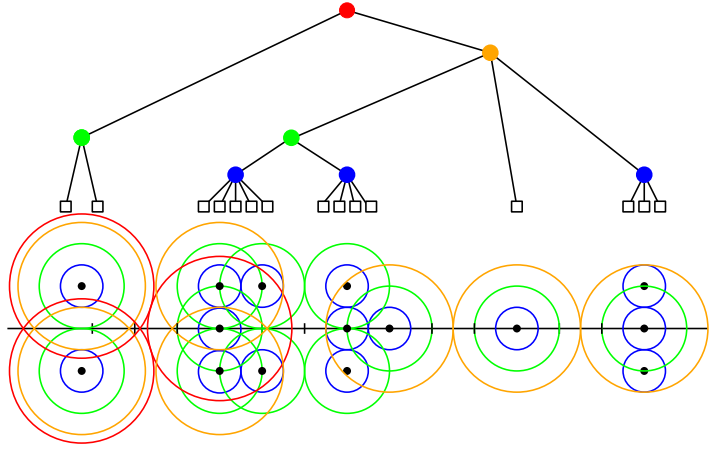
\vfigpdf{A dense pointset $P$ and construction of $\caltp$. Circles of different colors correspond
to different $\sigma$'s. The first choice of $\sigma$ corresponds to blue colored circles, 
followed by green, orange and red. 
The components formed are shown in the corresponding colors.
We only draw some of the relevant circles to give an idea.}{ripples}{0.5}

Let $u, v \in \caltp$ be such that $v$ is a child of $u$. We have the following properties of $\caltp$:
% (P1) There are exactly $\ell$ levels, and every level has a node with more than one children,
% because there are pair of roots with separation exactly $\sigma_i$.\\
% (P2) The component associated with the root node is $\callc$, because the maximum separation
% is $\sigma_\ell$. \\
% (P3) $\Delta_v \ge \sigma_u$. 
\begin{tightenum}{P}
\item$\sigma_u \le \min_{p\in G_v, q\in P\sm G_v}|p-q| \le 3 r_v$. The upper bound follows from the density of $P$.
The lower bound follows from the observation that the discs with radius
$\sigma/2$ centered at $p \in G_v$, where $\sigma_v <\sigma < \sigma_u$, do not touch
the discs of any other component, except when the radius is $\sigma_u/2$.
\item $r_u \le |G_u|\sigma_u$. Consider the graph $\calG$ with the vertices as $G_u$ and edges
between two vertices $p, q$ if $D(p, \sigma_u/2)\si D(q, \sigma_u/2)\neq \es$.
As $G_u$ is a connected component of these discs, we know that $\calG$ is connected.
Therefore, if $m$ is the number of vertices on the path joining $p, q$ in $\calG$,
then by triangular inequality $|p-q| \le m\sigma_u\le |G_u|\sigma_u$.
\item  If $p$ is a leaf-child of $u$ then $\sigma_u = \sigma_p$.
It is clear that any disc $D(p,r)$, with $r< \sigma_p/2$, cannot touch $D(q,r)$, for any other point $q$.
The first time they touch is when $\sigma_u=\sigma_p$.
If $p \in \CC \sm \RR$, then we further obtain that $|\Im(p)| \ge \sigma_p \ge \sigma_u$.
\item  The size of $\caltp= O(|P|)$. % There are exactly $\ell$ levels, and
Every level has a node with more than one child, as 
there are pairs of components with separation exactly $\sigma$.
\item  $P$ is the component associated with the root of $\caltp$.
\end{tightenum}

The next result is an amortization analogous to that of the Davenport-Mahler bound over the
root separation bound.
\bleml{dense} 
If $P$ is a dense pointset then
        \beql{mcint}
        \int_{(2D_P \si \RR) \sm \su_{p\in P} J_p} \frac{dx}{d(x,P)} = O(|P| \ln |P|),
        \eeql
where for $p \in P \si \RR$, $J_p \as [p\pm \sigma_p/2]$, and $J_p=\es$ otherwise.
\eleml
\bpf
We break the integral recursively over the nodes of $\caltp$. 
For an internal node $u$  of $\caltp$, we will show the following claim:
        \beql{intu}
        \int_{(2D(G_u)\si \RR )\sm \su_{p \in G_u}J_p} \frac{dx}{d(x,P)} = O(\nu_u \ln |G_u|).
        \eeql
We take the sum over all internal nodes $u$. From (P4) we know that
$|\caltp|=O(|P|)$, and hence $\sum_u \nu_u=O(|P|)$; moreover, 
from (P5) we know that the component associated with the root of $\caltp$ is $P$.
These observations  then give us \refeq{mcint}.

For a point $p \in P$, recall that $\calv_p$ is the set of points in $2D_P\si \RR$ 
closer to $p$ than to any other point of $P$; by definition $J_p \ib \calv_p$.  
Suppose $u$ is the parent of $p$.
We will bound the integral over $\calv_p$ in two steps: 
the portion of $\calV_p$ inside $2D(G_u)$ 
and the portion outside $2D(G_u)$. The latter portion is where amortization occurs, as 
for an $x \nin 2D(G_u)$, the distance of $x$ to $p \in G_u$ is roughly $|x- m_u|$.
Let $v$ be a child of $u$. There are three cases to consider:
\begin{tightenum}{}
\item {\bf Case 1.} $v$ is a leaf $p\in \RR$. 
We first bound the portion $I_p$ of $\calv_p$ inside $2D(G_u)$;
the portion outside will be handled collectively for all points in the third case.
% Let $J_p \as D(p, \sigma_p/2) \si \RR$. 
% Partition the interval $I_p = J_p \su (I_p \sm J_p)$;
% note that the interval $J_p$ always belongs to the Voronoi region of $p$. Correspondingly,
% we have
%         \beql{int1}
%         \int_{I_p}\gxdx= \int_{J_p}\gxdx+ \int_{I_p \sm J_p}\gxdx.
%         \eeql
% Consider the first of these integrals. % If $p$ is a midpoint of a maximal ssc then we have already
% % derived a bound on this in \refeq{bound1}. So suppose $p$ is a real root.
% For $x \in J_p$, our stopping function $G(x)=2/d_2(x,P)$.
% Suppose $q \in P$ is such that $d_2(p, P)=|p-q|$.
% Then for all $x \in J_p$, $d_2(x,P) \ge |p-q|-|p-x|\ge \sigma_p/2$,
% and hence we have the following upper bound 
%         $\int_{J_p} \frac{2dx}{d_2(x,P)} \le 4\int_{p-\sigma_p/2}^{p+\sigma_p/2} \frac{dx}{\sigma_p}=4$.
% Now consider the second integral on the RHS of \refeq{int1}.
For all $x \in I_p \sm J_p$, it is clear that $d(x, P) = |x-p|$. From
\refLem{int}(Re) we obtain that
        $\int_{I_p \sm J_p} \frac{dx}{|x-p|}=O\paren{\ln \frac{w(I_p)}{\sigma_p}}$.
But as $I_p \ib 2D(G_u) \si \RR$, we know that $w(I_p) \le 4r_u$.
From (P3), we know that $\sigma_p = \sigma_u$.
Therefore, $\int_{I_p \sm J_p} \frac{dx}{|x-p|} = O(\ln r_u/\sigma_u) = O(\ln |G_u|)$, from (P2).
% Hence $\int_{I_p\sm J_p} \frac{dx}{d(x,P)}= O(\ln |G_u|)$.
\item {\bf Case 2.} $v$ is a leaf $p=\Re(p)+i\Im(p) \in \CC \sm \RR$. Again consider the interval
$I_p \as \calv_p \si 2D(G_u)$; in this case $J_p=\es$. As $p$ is the closest point to any  $x \in I_p$,
$d(x, P)=|x-p|$. Moreover, $p$ and both the endpoints of $I_p$ are
in $2D(G_u)$, so the maximum distance of an endpoint of $I_p$ from $p$ is $\le 2r_u$.
Therefore, from \refLem{int}(Im) we have
        $$\int_{I_p} \frac{dx}{d(x,P)} = O\paren{\ln \frac{r_u}{|\Im(p)|}}.$$
But recall from (P3) that $|\Im(p)| \ge \sigma_u$, hence $r_u/|\Im(p)| \le r_u /\sigma_u \le |G_u|$,
where the last inequality follows from (P2). Therefore, $\int_{I_p} \frac{dx}{d(x,P)} = O(\ln |G_u|)$.
\item {\bf Case 3.} $v$ is an internal node. Inductively, we have already bounded
the integral $\int_{(2D(G_v)\si \RR) \sm \su_{p\in G_v}J_p}dx/d(x,P)$. However, it is possible that $\calV_p$, for
some point $p\in G_v$ extends beyond $2D(G_v)\si \RR$. 
Suppose $p$ is such a point, and $x\in W_p \as  \calV_p \si (2 D(G_u) \sm 2D(G_v))$. 
Then we know that $|x- p| \ge |x- m_v|/2$, where $m_v$ is the center of $D(G_v)$.
Therefore, 
        \begin{align*}
        \sum_{p \in G_v}\int_{W_p}\frac{dx}{|x-p|}   \le \int_{(2D(G_u) \sm 2D(G_v))\si \RR}\frac{2dx}{|x-m_v|}.    \end{align*}
As $2w(I_u) = 4r_u$, from \refLem{int}(Re), it follows that the integral on the RHS is bounded 
by $O(\ln r_u/r_v)$. But from (P2) we have $r_u \le |G_u|\sigma_u$, and $\sigma_u\le 3r_v$ from (P1).
Therefore, we obtain
\begin{align*}\sum_{p \in G_v}\int_{W_p}\frac{dx}{d(x,P)}=O(|\ln |G_u|).\end{align*}
This is the case where the amortization of the integral over the Voronoi regions takes place.
\end{tightenum}

Summing the bounds for all children $v$ of $u$ gives \refeq{intu}.
\epf

The following is the analogue of \refLem{dense} in  $\CC$:
define $D_p \as D(p, \sigma_p/2)$, then
        $$\int_{2D_P \sm \su_p D_p }\frac{dz}{d(z,P)} = O(|P|\ln |P|).$$

\ignore{
\subsection{Bound for Newton+EVAL}\label{sec:neval}
In this section, we derive an upper bound on the integral in the RHS of \refeq{nie3}
where the stopping function corresponds to the centered form interval arithmetic
based predicates used in the EVAL algorithm \cite{burr:contamort:13,sharma-yap:near-optimal:12}:
        \beql{evalp}
        \begin{split}
        C_0(I) \equiv |f(m(I))| > \sum_{j \ge 1}\abs{\frac{f^{(j)}(m(I))}{j!}} \paren{\frac{w(I)}{2}}^j\text{, and}\\          
        C_1(I) \equiv |f'(m(I))| > \sum_{j \ge 1}\abs{\frac{f^{(j+1)}(m(I))}{j!}} \paren{\frac{w(I)}{2}}^{j-1}.
        \end{split}
        \eeql
The stopping function in this case is $G(x) \as 1.5\min\set{S_0(x), S_1(x)}$, where
        \beql{seval}
        S_0(x) \as \sum_{\alpha \in Z(f)} \frac{1}{|x-\alpha|}
        \text{, and }
        S_1(x) \as \sum_{\alpha' \in Z(f')} \frac{1}{|x-\alpha'|}.
        \eeql
The idea is similar to what was done in \refSec{ndsc}, namely to charge $S_0$ on the region between
clusters, and $S_1$ on the roots. However, there is an added complication as $S_0$ depends on all
the roots, and not just the nearest root; similarly, for $S_1$ which depends on all critical points.
The bound is roughly $n$ times the bound in \refThm{bound}. This is because earlier only the
distance to a nearest root to $x$ played a significant role for $C_0$; however, $S_0$ is governed
by the distance to all the roots, and if all of them are equidistant from 
$x$ then $S_0 \sim n/d(x,V)$, which gives us the additional factor of $n$. In fact, the analysis
reveals that such a equi-distributed geometry of roots can perhaps achieve the bound.
We show the following:
\bthml{neval}
Given a square-free polynomial $f\in \RR[x]$ of degree $n$, 
the size of the subdivision tree constructed by \texttt{Newton-Isol}($f, I_0$) using
predicates given in \refeq{evalp} is bounded by $O(n^2\ln n)$.
\ethml

\bpf

\epf
}

\section{Concluding Remarks}
Our aim has been to devise a general approach to 
improve any subdivision based algorithm for real root isolation.
This is achieved by the \nie predicate, which detects strongly separated clusters,
and hence reduces the number of subdivisions from $O(\log R_\calc/r_\calc)$ 
to $O(n\log n)$. The crucial ingredient is Ostrowski's criterion based on
deviations of the Newton diagram of a polynomial. The criterion works
for complex polynomials, so we expect an analogue of \nisol for isolating
complex roots that is conceptually simpler than the existing approaches.
We have not explored the practical aspects of the algorithm, nevertheless,
we think that the analysis based on the geometry of cluster provides tools and
techniques for an alternate approach to understand existing algorithms.

%We derived bounds on the size of the subdivision tree of the algorithm.
We can bound the arithmetic complexity of \nisol as follows.
The Newton diagram computation takes $O(n)$, and the Taylor shift
$O(n \log n)$ operations. The number of Newton iterations
to approximate $\calc$ is bounded by $O(\log \log \frac{R_\calc}{r_\calc})$,
which is $O(\log (nL))$ from root separation bounds.  
Therefore, the arithmetic complexity, ignoring poly-log factors, 
is bounded by $\wt{O}(n^2)$. The extension
to the bitstream model involves deriving a robust
version of Ostrowski's result and bounding precision requirements.
The latter will be governed by perturbation bounds for clusters.
For a cluster of size $k$, we expect these bounds to be
$O(\eps^{1/k})$, for $\eps$-perturbation in the coefficients.
In the worst case, this would give an $O(n(L+\log n))$ bound on 
the precision.
%\eRemark

% \bibliographystyle{alpha}
% \bibliography{st,yap,exact,geo,alge,math,com,rob,cad,algo,visual,gis,quantum,mesh,tnt,fluid,kamath} % test.bib is a default file in MiKtex/bibtex/bib, used for testing

\newpage
\section*{Appendix}
We give the proof of \refLem{correctness}; the arguments are based standard manipulations
with Taylor series in alpha-theory of Smale et al. We first prove \refLem{correctness}(i), for
% Given $z \in \CC$, define $f_j(z) \as f^{(j)}(z)/j!$, i.e., the coefficient of $x^j$
% in $f(x+z)$. From \refeq{rkr}
%         \beql{rkrk}
%         \rho_k(z) = \max_{j<k} \abs{\frac{f_j(z)}{f_k(z)}}^{\frac{1}{(k-j)}},
%         \rho_{k+1}(z)= \min_{j>k} \abs{\frac{f_k(z)}{f_j(z)}}^{\frac{1}{(j-k)}}.
%         \eeql
which we need the following functions from \cite{bcss:bk} defined for $\fk$:
        \beql{gammakz}
        \beta_k(z)= \abs{\frac{\fk(z)}{f^{(k)}(z)}},\;
        \gamma_k(z)= \max_{j \ge 1}\abs{\frac{f^{(k+j)}(z)}{(j+1)! f^{(k)}(z)}}^{\frac{1}{j}}
        \eeql
and $\alpha_k(z) =\beta_k(z)\gamma_k(z)$.
We derive relations between these quantities and $\rho_k(z)$'s given in \refeq{rkr}.
Considering the RHS of $\rho_k(z)$ for $j=k-1$, we immediately have
        \beql{betakbound}
        \rho_k(z) \ge \abs{\frac{f_{k-1}(z)}{f_k(z)}}
                          = k \beta_k(z).
        \eeql
Multiplying and dividing the inner term on the RHS of $\gamma_k(z)$ in \refeq{gammakz}
by $(k+j)!/k!$ we obtain that
\begin{align*}
        \gamma_k(z) %&= \max_{j \ge 1} \paren{\frac{(k+j)!}{k!(j+1)!}}^{1/j}\abs{\frac{f_{k+j}(z)}{f_{k}(z)}}^{1/j}\\
                                 &\le \max_{j \ge 1} \paren{\frac{(k+j)!}{k!(j+1)!}}^{1/j}\max_{j > k}\abs{\frac{f_{j}(z)}{f_{k}(z)}}^{1/(j-k)}\\
                                 &= \max_{j \ge 1} \paren{\frac{(k+j)!}{k!(j+1)!}}^{1/j} \frac{1}{\rho_{k+1}(z)}.  
\end{align*}
% We now bound the maximum value that the supremum on the RHS can take. The fraction
% $$\frac{(k+j)!}{k!(j+1)!} = \frac{1}{j+1} {k+j \choose j} = \frac{1}{j+1} \frac{(k+j)(k+j-1) \ldots (k+1)}{j!} = \frac{1}{j+1}\prod_{\ell=1}^{j} \paren{\frac{k+\ell}{\ell}} 
%                                       \le \frac{(k+1)^j}{j+1}.
% $$
The max-term is bounded by $(k+1)$,
%        $$\max_{j \ge 1} \paren{\frac{(k+j)!}{k!(j+1)!}}^{1/j} \le k+1,$$
which implies that 
        \beql{gammakbound}
        \gamma_k( z) \rho_{k+1}(z) \le (k+1).
        \eeql
% By similar manipulations, we can show that $\rho_{k+1}(z) \gamma_k(z) \ge 1$.
Multiplying \refeq{betakbound} and \refeq{gammakbound} we obtain 
that  $\alpha_k(z) \Delta_k(z) \le 2.$ Therefore, if $\Delta_k(z) \ge 12$ then
$z$ is an approximate zero of $\fk$ with associated root in $D(z, 1.5\beta_k(z))\ib D(z,\frac{3\rho_k(z)}{2k})$,
where the inclusion follows from \refeq{betakbound}. The claim on
Newton iterates follows from \cite[p.~160,Thm.2]{bcss:bk}.
% $\alpha_k(z) \Delta_k(z) \le (k+1)/k \le 2$, as in 
% We will also need a lower bound. Again \footnote{These inequalities are just specialization of Proposition 1.7(c)
% of Giusti et al. (specifically, substitute $\ell=m$)}.
% We will often need the following result \cite[p.~161, Lem.~3]{bcss:bk}:

We now prove \refLem{correctness}(ii). We will need the following result \cite[p.~161, Lem.~3]{bcss:bk}: 
for a $u \in [0, 1)$
        \beql{binom}
        \sum_{j \ge 0} {k+j \choose j}u^j = \frac{1}{(1-u)^{k+1}}.
        \eeql
Let $\eps \as 1.5$, $\delta \as \eps\rho_k/k$, and $u \as \frac{\delta}{\rho_{k+1}}=\frac{\eps}{k\Delta_k}$;
here we express $\rho_k(z)$ by $\rho_k$ (similarly, for the other quantities). 
\bleml{fklb}
If $z$ is such that $\Delta_k(z) \ge 16$ then for $z' \in D(z, \delta)$, we have 
% $u \as \delta/\rho_{k+1}(z) < 1$ then $z' \in D(z, \delta)$. 
%        \beql{fk}
        $|f_k(z')| \ge |f_k(z)|(1-u)^{-(k+1)}/2$.
%        \eeql
%where $\phi_k(x) \as 2(1-x)^{k+1}-1$.
\eleml
\bpf
Take the absolute values in the Taylor expansion of $f^{(k)}(z')$ % we obtain that
        % $$f_k(z') = \frac{1}{k!}\paren{\sum_{j \ge 0} \frac{f^{(k+j)}(z)}{j!} (\pm\delta)^j}.$$
and apply the triangular inequality to obtain
        $$\abs{f_k(z')} \ge \frac{1}{k!} \paren{|f^{(k)}(z)| - \sum_{j \ge 1} \abs{\frac{f^{(k+j)}(z)}{j!}} \delta^j}.$$
Dividing both sides by $|f_k(z)|$,  and
multiplying and dividing the summation term on the RHS by $k!$ and $(k+j)!$ 
we obtain that 
%         $$\abs{\frac{f_k(z')}{f_k(z)}} \ge \paren{1 - \sum_{j \ge 1} \abs{\frac{f^{(k+j)}(z)}{j!f^{(k)}(z)}} \delta^j}.$$
% Multiplying and dividing the summation term on the RHS by $k!$ and $(k+j)!$ we obtain that
        $$\abs{\frac{f_k(z')}{f_k(z)}} \ge  \paren{1 - \sum_{j \ge 1} {k+j \choose j}\abs{\frac{f_{k+j}(z)}{f_{k}(z)}} \delta^j}.$$
From the expression of $\rho_{k+1}(z)$ in \refeq{rkr} and definition of $u$, we deduce that 
\begin{align*}
        \abs{\frac{f_k(z')}{f_k(z)}}
        %               &\ge |f_k(z)| \paren{1 - \sum_{j \ge 1} {k+j \choose j}\gamma_k(z)^j \delta^j}\\
        %               &= |f_k(z)| \paren{1 - \sum_{j \ge 1} {k+j \choose j}u^j}\\
                            &\ge \paren{2 - \sum_{j \ge 0} {k+j \choose j}u^j}.  
\end{align*}
%Since $u < 1/k$, we use \refeq{binom} to simplify the RHS to
%        $$\abs{f_k(z')} \ge |f_k(z)| \paren{2 - \frac{1}{(1-1/k)^{k+1}}} = |f_k(z)| \frac{\phi_k(1/k)}{(1-1/k)^{k+1}}.$$
Using \refeq{binom}, and the bound on $\Delta_k$, the RHS can be simplified to $(1-u)^{-(k+1)}/2$.
%        $$\abs{\frac{f_k(z')}{f_k(z)}}\ge \frac{\phi_k(u)}{(1-u)^{k+1}}.$$
\epf

\blem
If $\Delta_k(z)\ge 16$ then for all $z' \in D(z, \delta)$
we have $\rho_k(z') < 2e^6 \rho_k(z)$ and $\rho_{k+1}(z') \ge \rho_{k+1}(z)/3$.
Therefore, $\Delta_k(z') \ge \Delta_k(z)/(6e^6)$.
\elem
\bpf
%We next derive upper bounds on $f_j(z')$
% Let 
% $\eps\as 3/2$, $\delta \as \eps \beta_k(z)$, $u \as \delta/\rho_{k+1}$.
% for $z' = z\pm \delta$, where $\delta \as \eps \beta(\fk,z)$
% and $\eps = 3/2$.  There are two cases to consider,
% namely $j < k$ and $j > k$. Let us consider the first case:
For $j < k$, take absolute values in the Taylor expansion of $f_j(z')$, apply triangular inequality,
and split the summation up to $k$ and beyond $k$, to get 
% that
%         $$f_j(z') = \frac{1}{j!}\paren{f^{(j)}(z) + \sum_{\ell =1}^{k-j} \frac{f^{(\ell+j)}(z)}{\ell!} \delta^\ell + \sum_{\ell > k-j} \frac{f^{(\ell+j)}(z)}{\ell!} (\pm\delta)^\ell }.$$
% Taking absolute values and applying triangular inequality it follows that
%\begin{multline}
\begin{align*}
      |f_j(z')| %&\le |f_{j}(z)| + \sum_{\ell =1}^{k-j} \abs{\frac{f^{(\ell+j)}(z)}{j!\ell!}} \delta^\ell + \sum_{\ell > k-j} \abs{\frac{f^{(\ell+j)}(z)}{j!\ell!}} \delta^\ell \\
                     = |f_{j}(z)| + \sum_{\ell =1}^{k-j} {\ell+j \choose \ell}\abs{{f_{\ell+j}(z)}} \delta^\ell + 
                                \sum_{\ell > k-j} {\ell+j \choose \ell}\abs{f_{\ell+j}(z)} \delta^\ell .  
\end{align*}
%\end{multline}
        % |f_j(z')| %&\le |f_{j}(z)| + \sum_{\ell =1}^{k-j} \abs{\frac{f^{(\ell+j)}(z)}{j!\ell!}} \delta^\ell + \sum_{\ell > k-j} \abs{\frac{f^{(\ell+j)}(z)}{j!\ell!}} \delta^\ell \\
        %              = |f_{j}(z)| + \sum_{\ell =1}^{k-j} {\ell+j \choose \ell}\abs{{f_{\ell+j}(z)}} \delta^\ell + \sum_{\ell > k-j} {\ell+j \choose \ell}\abs{f_{\ell+j}(z)} \delta^\ell .  
Divide by $|f_k(z)|$ and use the expressions in \refeq{rkr} to obtain
        $$\abs{\frac{f_j(z')}{f_k(z)}}\le \rho_k^{k-j}+ \sum_{\ell =1}^{k-j} {\ell+j \choose \ell}\rho_k^{k-\ell-j} \delta^\ell + \sum_{\ell > k-j} {\ell+j \choose \ell}\frac{\delta^\ell}{\rho_{k+1}^{\ell+j-k}}.$$
Since $\delta =\eps\rho_k/k$, we can pull 
out $\rho_k^{k-j}$ from the RHS (and since $\Delta_k>1$), we get that
%\begin{multline}
$$
        \abs{\frac{f_j(z')}{f_k(z)}}%&\le \beta_k^{k-j}\paren{1+ \sum_{\ell =1}^{k-j} {\ell+j \choose \ell}\paren{\frac{\eps}{k}}^{\ell} + \frac{1}{\alpha_k^{k-j}}\sum_{\ell > k-j} {\ell+j \choose \ell}u^\ell}\\
                                               %&\le \beta_k^{k-j}\paren{1+ \sum_{\ell =1}^{k-j} {\ell+j \choose \ell}\paren{\frac{\eps}{k}}^\ell + \frac{1}{\alpha_k^{k-j}}\sum_{\ell > k-j} {\ell+j \choose \ell}\paren{\frac{\eps\alpha_k}{k}}^\ell}\\
                                               %&\le \rho_k^{k-j}\paren{1+ \sum_{\ell =1}^{k-j} {\ell+j \choose \ell}\paren{\frac{\eps}{k}}^\ell  + \sum_{\ell > k-j} {\ell+j \choose \ell}\Delta_k^{-(\ell-k+j)}\paren{\frac{\eps}{k}}^\ell }\\
                                               \le \rho_k^{k-j}\paren{1+ \sum_{\ell =1}^{k-j} {\ell+j \choose \ell}\paren{\frac{\eps}{k}}^\ell  
                                                        + \sum_{\ell > k-j} {\ell+j \choose \ell}\paren{\frac{\eps}{k}}^\ell}.
%                                               &=  \paren{\frac{\beta_k}{\alpha_k}}^{k-j}\paren{\alpha_k^{k-j}+ \sum_{\ell =1}^{k-j} {\ell+j \choose \ell}\alpha_k^{k-j} k^{-\ell} +\sum_{\ell > k-j} {\ell+j \choose \ell}u^\ell}.
$$ %\end{multline}
% where in the last step we use the fact that $\alpha_k \le 1$.
Assuming $k \ge 2$, from \refeq{binom} we obtain that
        $$        \abs{\frac{f_j(z')}{f_k(z)}}
                                               \le \rho_k^{k-j}\paren{1- \frac{\eps}{k}}^{-(j+1)}.$$
% Since $\alpha_k < 1$, the term $\alpha_k^{k-j}$ in the second summation term on the RHS can be
% replaced by $\alpha_k^{\ell}$ to give
% $$        \abs{\frac{f_j(z')}{f_k(z)}}
%                                                \le  \paren{\frac{\beta_k}{\alpha_k}}^{k-j}\paren{\alpha_k^{k-j}+ \sum_{\ell =1}^{k-j} {\ell+j \choose \ell}\paren{\frac{\alpha_k}{k}}^{\ell} +\sum_{\ell > k-j} {\ell+j \choose \ell}u^\ell}.$$
% But $u=\gamma_k\delta \le \alpha_k/k$, and hence we obtain that
%         $$\abs{\frac{f_j(z')}{f_k(z)}}\le  \paren{\frac{\beta_k}{\alpha_k}}^{k-j}\paren{\alpha_k^{k-j}+ \sum_{\ell \ge 1} {\ell+j \choose \ell}\paren{\frac{\alpha_k}{k}}^{\ell}}.$$
% Applying \refeq{binom} to the RHS we get
%         $$\abs{\frac{f_j(z')}{f_k(z)}}
%                 \le  \paren{\frac{\beta_k}{\alpha_k}}^{k-j}\paren{\alpha_k^{k-j}+ \paren{\frac{1}{1-\alpha_k/k}}^{k+1}-1}
%                 \le  \paren{\frac{\beta_k}{\alpha_k}}^{k-j}\paren{\frac{1}{1-\alpha_k/k}}^{k+1},$$
% as $\alpha_k < 1$. 
Combining this bound with \refLem{fklb}, and doing some further simplifications
%and the observation that
%as $u$ increases to $\eps\alpha_k/k$ the ratio $\phi_k(u)/(1-u)^{k+1}$ decreases, 
we obtain the upper bound on $\rho_k(z')$. Note that we require $k \ge 2 > \eps$.
%         $$\abs{\frac{f_j(z')}{f_k(z')}}^{1/(k-j)} %= \paren{\abs{\frac{f_j(z')}{f_k(z)}} \abs{\frac{f_k(z)}{f_k(z')}}}^{1/(k-j)}
%                 \le \rho_k \paren{\frac{(1-\eps/(k\Delta_k))^{k+1}}{(1-\eps/k)^{j+1}}}^{1/(k-j)} \paren{\frac{1}{\phi_k(\eps/(k\Delta_k))}}^{1/(k-j)}.$$
% Since $\phi_k(\eps/(\Delta_kk)) < 1$, the second term in parenthesis on the RHS is smaller than 
% the inverse of $\phi_k(\eps/k\Delta_k)$. The quantity on the RHS can be shown to be smaller than $2e^6$.
To derive a lower bound on $\rho_{k+1}(z')$ in terms of $\rho_{k+1}(z)$,
% we need an upper bound on $|f_j(z')|$, where $j > k$. 
we take absolute values in the Taylor expansion of $f_j(z')$, for $j >k$, 
apply the triangular inequality,  and divide both sides by $|f_k(z)|$, to get
        $$\abs{\frac{f_j(z')}{f_k(z)}} %\le \abs{\frac{f_j(z)}{f_k(z)}} + \sum_{\ell \ge 1} \abs{\frac{f^{(\ell+j)}(z)}{j!\ell!f_k(z)}} \delta^\ell
                                                      \le \abs{\frac{f_j(z)}{f_k(z)}} + \sum_{\ell \ge 1} {\ell+j\choose \ell}\abs{\frac{f_{\ell+j}(z)}{f_k(z)}} \delta^\ell.$$
From the expression for $\rho_{k+1}$ in \refeq{rkr} and \refeq{binom} it follows that
        $$\abs{\frac{f_j(z')}{f_k(z)}} %\le \gamma_k^{j-k} + \sum_{\ell \ge 1} {\ell+j\choose \ell}\gamma_k^{\ell+j-k}\delta^\ell 
                                                     % = \gamma_k^{j-k}\paren{1 + \sum_{\ell \ge 1} {\ell+j\choose \ell}\paren{\frac{\delta}{\rho_{k+1}}}^\ell}
                                                      \le \frac{1}{\rho_{k+1}^{j-k}} \paren{1-u}^{-(j+1)}.$$
% where in the last step we use \refeq{binom}. 
% Since
%         $$\abs{\frac{f_j(z')}{f_k(z')}}^{1/(j-k)} = \paren{\abs{\frac{f_j(z')}{f_k(z)}} \abs{\frac{f_k(z)}{f_k(z')}}}^{1/(j-k)}$$
Combining this with the lower bound in \refLem{fklb}, and using the lower bound on 
$\Delta_k$ we further obtain that
        $$\abs{\frac{f_j(z')}{f_k(z')}}^{1/(j-k)} \le \frac{2^{j-k}}{(1-u)\rho_{k+1}}.$$% \paren{\frac{1}{\phi_k(u)}}^{1/(j-k)}.$$
Since $u \le \eps/(k\Delta_k)$ and $\Delta_k \ge16$, we get that
% $\phi_k(u)> 1/2$ and hence
        $$\abs{\frac{f_j(z')}{f_k(z')}}^{1/(j-k)} \le \frac{3}{\rho_{k+1}},$$
which implies the desired lower bound on $\rho_{k+1}(z')$. % \ge \rho_{k+1}(z)/3$, for all $z' \in \DD(z, \eps\beta(\fk, z))$.
\epf

To show \refLem{correctness}(iii), we suppose that $\rho_k(w) \le \rho_k(z)$.
As the two inclusion discs intersect it follows that
        $$|w \pm 3\rho_k(w) -z| \le |z-w|+3\rho_k(w) \le 9\rho_k(z) \le \frac{\rho_{k+1}(z)}{3},$$
where the last inequality follows from $\Delta_k(z)\ge 27$. This implies that
$D(w, 3\rho_k(w)) \ib D(z, \frac{\rho_{k+1}(z)}{3})$, and hence both the discs have the same cluster.

We next give a self-contained proof of Pawlowski's result \cite[Thm.~2.2]{pawlowski:zeros-of-derivatives:99}.
The difference in our proof is that we avoid using Enestr\"om-Kakeya theorem. The key idea of using
Walsh's representation theorem \cite[Thm.~3.4.1c]{rahman-schmeisser:polynomials:bk}, however, is common to both the proofs.

Given a cluster $\calc$ of size $k$, and $z\in \CC$,
let $f_1(z)=\lead(f)\prod_{\alpha \in \calc}(z-\alpha)$ and $f_2(z) =\prod_{\beta \in \calcc}(z-\beta)$.
From Leibniz's formula we have
        $$f^{(j)}(z)=\sum_{i=\max\set{0, j+k-n}}^{\min\set{j,k}} {j \choose i} f_1^{(i)}(z) f_2^{(j-i)}(z).$$
We will focus on the case when $j \le k$, in which case the upper bound of the summation is $j$.
The bounds on the summation are required because $f_1$ cannot be differentiated more than $k$ times and
similarly for $f_2$. % We will mostly be interested in the case when $j \le k$, which simplifies the equation above
% to
%         $$f^{(j)}(z)=\sum_{i=\max\set{0, j+k-n}}^{j} {j \choose i} f_1^{(i)}(z) f_2^{(j-i)}(z).$$
Applying Walsh's representation theorem, first to $f_1$ we obtain that there is an $\alpha \in D(m,r_\calc)$ such that
        $$f^{(j)}(z)=\sum_{i=\max\set{0, j+k-n}}^{j} {j \choose i} k(k-1) \cdots (k-i+1) (z-\alpha)^{k-i}f_2^{(j-i)}(z).$$
Now applying Walsh's representation theorem to $f_2$, we know that there is a $\beta \nin D(m,R_\calc)$ such that
        $$f^{(j)}(z)=\sum_{i=\max\set{0, j+k-n}}^{j} {j \choose i} k(k-1) \cdots (k-i+1) (z-\alpha)^{k-i}
                                                                                        (n-k)(n-k-1)\cdots (n-k-j+i+1) (z-\beta)^{n-k-j+i}.$$
Opening up the binomial  ${j \choose i}$ and simplifying we obtain
        $$f_j(z) = \frac{f^{(j)}(z)}{j!}=\sum_{i=\max\set{0, j+k-n}}^{j} {k\choose i} (z-\alpha)^{k-i}
                                                                                        {n-k\choose j-i} (z-\beta)^{n-k-j+i}.$$
Pulling out the last term from the RHS we obtain that 
        $$f_j(z)={k \choose j}(z-\alpha)^{k-j}(z-\beta)^{n-k}\sum_{i=\max\set{0, j+k-n}}^{j} \frac{{k\choose i} {n-k\choose j-i}}{{k \choose j}}\paren{\frac{z-\alpha}{z-\beta}}^{j-i};$$
note that when $j=k$, the term $(z-\alpha)^{k-j}$ is one. Now we  substitute $i$ by $j-i$ to obtain
        \beql{fjz}
        f_j(z)={k \choose j}(z-\alpha)^{k-j}(z-\beta)^{n-k}\sum_{i=0}^{\min\set{j, n-k}} \frac{{k\choose j-i} {n-k\choose i}}{{k \choose j}}\paren{\frac{z-\alpha}{z-\beta}}^{i}.
        \eeql
% Our claim is that in a suitable disc around a strongly separated cluster the first term in the summation above 
% will dominate the rest in absolute value and hence $f_j(z)$ cannot vanish. 
The fraction
$$\frac{{k\choose j-i} {n-k\choose i}}{{k \choose j}} 
        = \frac{1}{i!}\frac{(n-k)!}{(n-k-i)!} \frac{j!}{(j-i)!} \frac{(k-j)!}{(k-j+i)!}
        \le \frac{1}{i!} \paren{\frac{(n-k)j}{(k-j)}}^i;$$
note that for $j=k$, the denominator does not appear; we capture this by using the notation
$(x)_1 \as \max\set{1, x}$.
% Clearly, the maximum value for this
% fraction is obtained at $j=k-1$, and it is smaller than the value at $j=k$,
% namely $(k(n-k))^i/i!$.
Therefore, the summation in \refeq{fjz}
does not vanish if $z$ satisfies the following inequality:
        $$\half\ge  \sum_{i\ge 1}^{\min\set{j,n-k}} \frac{1}{i!}\paren{\frac{j(n-k)}{(k-j)_1}\abs{\frac{z-\alpha}{z-\beta}}}^{i}.$$
Substituting the upper bound of the summation by infinity, we get the following stronger constraint:
        $$\half \ge  \sum_{i\ge 1} \frac{1}{i!}\paren{\frac{j(n-k)}{(k-j)_1}\abs{\frac{z-\alpha}{z-\beta}}}^{i}.$$
Adding one on both sides, and using the observation that the RHS is the expansion of 
$\exp(\paren{\frac{j(n-k)}{(k-j)_1}\abs{\frac{z-\alpha}{z-\beta}}})$, we get that the inequality follows
if
        \beql{ab}
        \ln 1.5 \ge \frac{1}{4} \ge \paren{\frac{j(n-k)}{(k-j)_1}\abs{\frac{z-\alpha}{z-\beta}}}.
        \eeql
% From the AM-GM inequality we have that $k(n-k)\le n^2/4$, therefore, the inequality above follows
% if
%         \beql{ab}
%         \ln 2> \paren{\frac{n^2}{4}\abs{\frac{z-\alpha}{z-\beta}}}.
%         \eeql
% We claim that the above inequality holds if $z$ is such that $|z-m_\calc|\le R_\calc/(2n^2)$.
% Let $M \as \max \set{|z-m_\calc|, r_\calc}$. 
Since $\alpha \in D(m_\calc, r_\calc)$, we have
        $$|z-\alpha| \le |z-m_\calc| + r_\calc \le 2\max(|z-m_\calc|,r_\calc),$$ 
and 
$$|z-\beta| \ge |m_\calc- \beta| - |z-m_\calc| \ge R_\calc - |z-m_\calc| \ge R_\calc - \max\set{|z-m_\calc|, r_\calc}.$$
%n^2\max(|z-m_\calc|,r_\calc)$
% (since $R_\calc \ge 2n^2 \max(|z-m_\calc|,r_\calc)$ for a strongly separated cluster).
These bounds imply that 
        $$\frac{|z-\alpha|}{|z-\beta|} \le 2 \frac{\max\set{|z-m_\calc|, r_\calc}}{(R_\calc - \max\set{|z-m_\calc|, r_\calc})},$$ 
and hence \refeq{ab} follows
if 
        $$R_\calc \ge \max\set{|z-m_\calc|, r_\calc} \paren{1 + \frac{8j(n-k)}{(k-j)_1}}.$$
To summarize, we have the following result:

\bleml{derivroots}
Let $\calc$ be a cluster of size $k$ and $j \in \set{0 \dd k}$.
If $z \in \CC$ is such that 
        \beql{rcalc}
        R_\calc \ge 2\max\set{|z-m_\calc|, r_\calc}\paren{\frac{8j(n-k)}{(k-j)_1}}_1,
        \eeql
%        $$|z-m_\calc| \le \frac{R_\calc}{2n^2}$$
then  %for $j \le k$,
        $$f_j(z) = {k \choose j}(z-\alpha)^{k-j}(z-\beta)^{n-k}(1 \pm \half),$$
for some $\alpha \in D(m, r_\calc)$ and $\beta \nin D(m,R_\calc)$;
the notation ``$\pm$''  stands for a $\theta\in \CC$ such that $|\theta| \le 1$.
Moreover, if $z$ also satisfies $|z-m_\calc| > r_\calc$ then $f^{(j)}(z) \neq 0$.
\elem

We specialize this result to the case of a strongly separated cluster:
\bcorl{derivrootsa}
For a strongly separated cluster $\calc$, if $z$ is such that
        $$|z-m_\calc| \le \frac{R_\calc}{ 4n^2}$$
then  for $0 \le j \le k$,
        $$f_j(z) = {k \choose j}(z-\alpha)^{k-j}(z-\beta)^{n-k}(1 \pm \half),$$
for some $\alpha \in D(m, r_\calc)$ and $\beta \nin D(m,R_\calc)$.
Moreover, if $z$ also satisfies $|z-m_\calc| > r_\calc$ then $f^{(j)}(z) \neq 0$, for $0 \le j \le k$.
\ecorl
\bpf
Note that the maximum value of the term $j(n-k)/(k-j)_1$ is obtained at $j=k$, and it is 
$k(n-k)$. From the AM-GM inequality, we know that $k(n-k) \le n^2/4$. Therefore,
\refeq{rcalc} follows if $z$ and $r_\calc$ are such that
        $$R_\calc \ge 4n^2\max\set{|z-m_\calc|, r_\calc}.$$
For a strongly separated cluster $\calc$ we know that $R_\calc \ge 4n^2r_\calc$, and
the condition on $z$ is the condition in the corollary.
\epf

We use this result to show that $(k-j)$ roots of the $j$th derivative are in $D(m_\calc,r_\calc)$ and the remaining are
outside $D(m_\calc, R_\calc/(2n^2))$. 
Let $\alpha_1 \dd \alpha_k$ be the roots of $f$ in $\calc$ and let $\beta_1 \dd \beta_{n-k}$
be the remaining roots. Let $g_t$ be the polynomial with roots 
$(1-t)m_\calc+t\alpha_1 \dd (1-t)m_\calc+t\alpha_k, \beta_1 \dd \beta_{n-k}$. Thus 
$g_0(z)=(z-m_\calc)^k\prod_{\beta \in \calcc}(z-\beta)$ and $g_1(z)=f(z)$.
Since $g_t(z)$ has a strongly separated cluster of size $k$ in $D(m_\calc, tr_\calc)$, from the lemma above we know that
$g_t^{(j)}(z)$ does not vanish on the boundary of $D(m_\calc, r_\calc)$. As the roots vary continuously
with $t$, and $g_0^{(j)}(z)$ has a root of multiplicity $(k-j)$ at $m_\calc$, it follows that $g_t^{(j)}(z)$ has
$k-j$ roots in $D(m_\calc, tr_\calc)$ and the remaining roots outside $D(m_\calc, R_\calc/(2n^2))$.
Substituting $t=1$ gives us the desired result. To summarize, we have obtained the following result:

\bleml{derivrootsb}
Given a strongly separated cluster $\calc$ of size $k$, for $j \le k$, there are $k-j$ roots of the derivative
$f^{(j)}(z)$ in $D(m_\calc, r_\calc)$ and the remaining $(n-k)$ roots are outside $D(m_\calc, R_\calc/(2n^2))$.
\eleml
% \Remark
% This result is an improvement over Renegar's and Coppersmith-Neff's result.
% Their results state  that the roots of $f^{(j)}$ close to the cluster are inside $\DD(m_\calc, nr_\calc)$.
% But this result is similar to the result of Pawlowski and one in Rahman-Schmeiser, where they derive
% inclusion discs for the roots of the derivative.
% \eRemark

\end{document}